\documentclass[review,hidelinks,onefignum,onetabnum]{siamart250211}

\usepackage{lipsum}
\usepackage{graphicx}
\usepackage{epstopdf}
\usepackage{amssymb,amsfonts,mathrsfs}
\usepackage{algorithm}
\usepackage{algpseudocode}

\ifpdf
  \DeclareGraphicsExtensions{.eps,.pdf,.png,.jpg}
\else
  \DeclareGraphicsExtensions{.eps}
\fi

\newsiamremark{remark}{Remark}
\newsiamremark{hypothesis}{Hypothesis}
\crefname{hypothesis}{Hypothesis}{Hypotheses}
\newsiamthm{claim}{Claim}
\newsiamremark{fact}{Fact}
\crefname{fact}{Fact}{Facts}

\headers{Efficient methods for the uncertain Boltzmann equation}{Y. Lin and L. Liu}

\title{Efficient numerical methods for the uncertain Boltzmann equation based on a hybrid solver\thanks{Submitted to the editors July 25, 2025.
\funding{The work of Y. Lin was partially supported by National Key R\&D Program of China (2020YFA0712000) and National Natural Science Foundation of China (12201404). L. Liu acknowledges the support by National Key R\&D Program of China (2021YFA1001200), Ministry of Science and Technology in China, Early Career Scheme (24301021) and General Research Fund (14303022 \& 14301423) funded by Research Grants Council of Hong Kong.}}}

\author{
Yiwen Lin\thanks{School of Mathematical Sciences, Shanghai Jiao Tong University, Shanghai 200240, P. R. China  (\email{linyiwen@sjtu.edu.cn}).}
\and 
Liu Liu\thanks{The Chinese University of Hong Kong, Hong Kong  (\email{lliu@math.cuhk.edu.hk}).}
}

\usepackage{amsopn}

\ifpdf
\hypersetup{
  pdftitle={Efficient methods for the uncertain Boltzmann equation},
  pdfauthor={Y. Lin and L. Liu}
}
\fi

\begin{document}

\maketitle

\begin{abstract}
In this work, we propose and compare several approaches to solve the Boltzmann equation with uncertain parameters, including multi-level Monte Carlo and multi-fidelity methods that employ an  asymptotic-preserving-hybrid (APH) scheme (Filbet and Rey, 2015) for the deterministic Boltzmann model. By constructing a hierarchy of models from finer to coarser meshes in phase space for the APH scheme and adopting variance reduction techniques, the MLMC method is able to allocate computational resources across different hierarchies quasi-optimally. On the other hand, in the bi-fidelity method we choose the APH scheme for the Boltzmann equation as the high-fidelity solver, and a finite volume scheme for the compressible Euler system as the low-fidelity model. Since both methods are non-intrusive, they can preserve the physical properties of the deterministic solver. Extensive numerical experiments demonstrate that our APH-based MLMC and multi-fidelity methods are significantly faster than standard approaches, while maintaining accuracy. 
We also provide practical guidelines for selection between APH-based MLMC and multi-fidelity approaches, based on solution smoothness and computational resource availability. 
\end{abstract}



\section{Introduction}
\label{sec1:intro}

The Boltzmann equation is one of the most fundamental equations in kinetic theory whose applications range from rarefied gas dynamics, plasma physics to semiconductor modeling, radiative transfer, social and biological system etc. It provides a rigorous statistical description of non-equilibrium gas dynamics that bridges microscopic particle interactions with macroscopic fluid behavior \cite{Cercignani1988}.
However, the complexity of numerical computation of the Boltzmann equation remains a significant challenge, particularly due to its high dimensionality, nonlinearity of the collision operator and multiscale nature.

A key dimensionless parameter governing the multiscale dynamics is the Knudsen number (Kn), defined as the ratio of the mean free path over a characteristic macroscopic length scale. In practical applications, the Knudsen number can vary by orders of magnitudes across the computational domain, making numerical computation prohibitively expensive due to the stiffness. Asymptotic-preserving (AP) schemes have been popularly used for hyperbolic type kinetic equations, which preserves the transition from kinetic to fluid regime at the discrete level \cite{Jin1999, Jin2022an}. 
When the Knudsen number becomes small (Kn $\rightarrow 0$), the numerical scheme automatically converges to a consistent discretization of the limiting fluid system (e.g., Euler or Navier-Stokes). Regarding the Boltzmann equation under the diffusive scaling, many AP methods have been studied, including the BGK-operator penalization \cite{Filbet2010jcp}, exponential integrator \cite{Dimarco2011sinum} and micro-macro decomposition \cite{Lemou2008sisc, Gamba2019jcp}.
On the other hand, there are substantial challenges in numerically evaluating the collision operator and solving high-dimensional problem in the phase space \cite{Pareschi2000sinum, Mouhot2006mc, Gamba2017sisc}. Recent work such as fast spectral methods \cite{Gamba2017sisc} and model reduction techniques \cite{Hu2024nsf} are developed to tackle these difficulties.

While AP methods provide uniform accuracy across different scales without resolving the discretizations, hybrid methods can be more efficient by dynamically coupling kinetic and fluid solvers in different spatial regions. The key challenge is designing a robust domain decomposition criterion to switch between models.
Several strategies have been proposed: 
the most common one is the local Knudsen number criterion \cite{Boyd1995pf, Kolobov2007jcp} for adaptive switching between computationally efficient (but limited) hydrodynamic models and more accurate kinetic solvers. This criterion has been recently further implemented \cite{Degond2012jcp}.
Besides, the hydrodynamic breakdown indicator \cite{Tiwari1998jcp} provides a similar criterion based on the viscous and heat fluxes of the Navier-Stokes equation. There are other works such as \cite{Degond2010jcp, Alaia2012jcp, Dimarco2014jcp, Dimarco2008mms}. 
Later the authors in \cite{Filbet2015sisc} introduced an improved hybrid scheme with automatic domain decomposition criterion that identifies accurately the fluid and kinetic zones. We will apply the methodology in \cite{Filbet2015sisc} for our deterministic solver, and in each kinetic cell, an AP scheme \cite{Filbet2010jcp} is adopted. 

In reality, uncertainties in kinetic models arise from empirical collision kernels, imprecise initial/boundary conditions and measurement errors in material parameters, thus it becomes a timely task to study their corresponding uncertainty quantification (UQ) problems. To resolve numerical challenges of UQ problems such as curse of dimensionality, multi-level Monte Carlo (MLMC) and multi-fidelity methods in the stochastic collocation framework \cite{Babuska2010siamrev, Narayan2012sisc} provide efficient ways. MLMC methods exploit hierarchical approximations to reduce variance efficiently \cite{Caflisch1998an, Giles2008or, Mishra2016siamjuq, Mishra2013springer}, whereas multi-fidelity approaches leverage low-fidelity models (e.g., Euler equations for small Knudsen numbers) to accelerate convergence \cite{Narayan2014sisc, Peherstorfer2018siamrev, Dimarco2019jcp, Dimarco2020mms, Liu2020jcp}.
Specifically, by decomposing into several hierarchical levels, MLMC method has shown to perform well for stochastic differential equations \cite{Abdulle2013jcp, Barth2011nm}, including kinetic systems \cite{Hu2023jcp}. 
On the other hand, thanks to less complex and computationally cheaper low-fidelity models for practical problems, numerous multi-fidelity approaches are developed by combining different models to efficiently approximate high-fidelity solutions, see  \cite{Zhu2014siamjuq, Peherstorfer2016cmame, Peherstorfer2016sisc, Eldred2017springer}, including recent applications to kinetic equations \cite{Dimarco2019jcp, Liu2020jcp} and the Schr\"odinger equation \cite{LinLiu2025sisc}. 

By incorporating the hybrid kinetic/fluid coupling strategies and AP scheme for the Boltzmann equation, the goal of this work is to develop an asymptotic-preserving-hybrid (APH) based, efficient numerical methods for solving the uncertain Boltzmann equation under multiple scalings, in particular we study and compare two non-intrusive approaches: multi-level Monte Carlo (MLMC) and multi-fidelity methods. 

For the MLMC approach, we construct a hierarchy of models from finer to coarser meshes in phase space by employing the APH scheme at each level, and uses variance reduction to allocate overall computational resources. On the other hand, in the bi-fidelity method, we employ the Euler equation as the low-fidelity model and APH scheme for the high-fidelity Boltzmann model, while a hierarchical structure consisting of the standard AP method, APH scheme for the Boltzmann, and Euler model (high-, medium- and low-fidelity) is introduced in our tri-fidelity approach. By conducting extensive numerical experiments, we demonstrate that these two methods can achieve significant improvements in computational efficiency while maintaining accuracy. The APH-based MLMC method delivers a prominent speedup compared to standard MC method and can accurately solve challenging discontinuous problems such as those containing shocks. Our multi-fidelity methods exhibit optimal trade-off between computational cost and accuracy for Boltzmann model under various regimes. Furthermore, we provide practical selection guidelines for these methods based on specific problems and computational resources, which hopefully will shed light on solving more general UQ problems for the multiscale kinetic equations. 


This paper is organized as follows. Section \ref{sec1:intro} introduces the application and challenges of the Boltzmann equation with uncertainties, along with the objectives of this work. Section 2 presents the uncertain Boltzmann equation and its hydrodynamic limit. Section 3 describes a hybrid numerical method for solving the Boltzmann equation with asymptotic preserving property. Section 4 discusses two types of stochastic collocation based methods for handling uncertainties in the Boltzmann equation. Section 5 validates the accuracy and efficiency of the APH solver, compares different types of non-intrusive uncertainty quantification approaches, and provides recommendations for uncertain Boltzmann equation simulations. Finally, Section 6 concludes the paper.

\section{Boltzmann equation}
\label{sec2:Boltzmann}

\subsection{The Boltzmann equation with uncertainties}

Let $(\Omega, \mathscr{F}, \mathbb{P})$ be a probability space with $\Omega$ being the set of elementary events, $\mathscr{F}$ the corresponding $\sigma$-algebra, and $\mathbb{P}$ the probability measure mapping $\Omega$ into $[0,1]$ such that $\mathbb{P}(\Omega)=$ 1. 

The dimensionless Boltzmann equation with uncertainties is modeled by the following
\begin{equation}\label{eq:B_sto}
\left\{\begin{array}{l}
\partial_t f^{\varepsilon}(\omega;\cdot)+v \cdot \nabla_x f^{\varepsilon}(\omega;\cdot)=\dfrac{1}{\varepsilon} \mathcal{Q}_{\mathcal{B}}(f^{\varepsilon}, f^{\varepsilon})(\omega;\cdot), \quad  \omega \in \Omega, x \in D, v \in \mathbb{R}^{d_v},\\
f^{\varepsilon}(\omega; 0, x, v)=f_0(\omega; x, v),
\end{array}\right.    
\end{equation}
where $\omega \in \Omega^{d_z}$ is a $d_z$-dimensional random parameter with probability distribution $\pi(\omega)$ known in priori characterizing the uncertainties in the system. 
$f^{\varepsilon}=f^{\varepsilon}(\omega; x, v, t)$ is the probability density distribution function, modeling the probability of finding a particle at time $t>0$, position $x \in D \subset \mathbb{R}^{d_x}$ and velocity $v \in \mathbb{R}^{d_v}$, where $d_x$ and $d_v$ are the dimensions of the $x$ and $v$ variables. The open set $D\subset \mathbb{R}^{d_x}$ is a bounded Lipschitz continuous domain and thus the model \eqref{eq:B_sto} has to be supplemented with boundary conditions described later. The parameter $\varepsilon$ is the Knudsen number defined as the ratio of the mean free path over a typical length scale.  The collision operator $\mathcal{Q}$ is a high-dimensional, quadratic integral operator modeling the binary elastic collision between particles, given by
$$
\begin{aligned}
    &\mathcal{Q}_{\mathcal{B}}(f, f)(\omega;v)\\
    &=\int_{\mathbb{R}^{d_v}} \int_{\mathbb{S}^{d_v-1}} B\left(\omega;\left|v-v_*\right|, \cos \theta\right)\left(f\left(\omega;v^{\prime}\right) f\left(\omega;v_*^{\prime}\right)-f(\omega;v) f\left(\omega;v_*\right)\right) d \sigma d v_*.
\end{aligned}
$$ 
Here, $\left(v, v_*\right)$ and $\left(v^{\prime}, v_*^{\prime}\right)$ are the velocity pairs before and after the collision, with momentum and energy conserved; thus ($v^{\prime}, v_*^{\prime}$) can be expressed in terms of ($v, v_*$) as follows:
$$
v^{\prime}=\frac{v+v_*}{2}+\frac{\left|v-v_*\right|}{2} \sigma, \quad v_*^{\prime}=\frac{v+v^*}{2}-\frac{\left|v-v_*\right|}{2} \sigma,
$$
with the vector $\sigma$ the scattering direction varying on the unit sphere $\mathbb{S}^{d_v-1}$.
The collision kernel $B\left(\omega;v, v_*, \sigma\right)=B\left(\omega;\left|v-v_*\right|, \cos \theta\right)$ is a non-negative function  which by physical arguments of invariance only depends on $\omega$, $\left|v-v_*\right|$ and $\cos \theta=\frac{\sigma \cdot\left(v-v_*\right)}{\left|v-v_*\right|}$. 
We consider the variable hard sphere (VHS) model \cite{Bird1994oxford}, with the collision kernel given by:$$B\left(\omega;v, v_*, \sigma\right)=b\left|v-v_*\right|^\gamma, \quad-d_v<\gamma \leq 1,$$where $b$ is a positive constant, $\gamma>0$ corresponds to the hard potential, and $\gamma<0$ is the soft potential.  

The collision operator $\mathcal{Q}$ has the fundamental properties of conserving mass, momentum and kinetic energy: 
\begin{equation}\label{eq:conservation}
\int_{\mathbb{R}^{d_v}} \mathcal{Q}(f, f) m(v) d v=0, \text{ with } m(v)=\left(1, v, \frac{|v|^2}{2}\right)^{\top},    
\end{equation}
and the well-known Boltzmann’s H-theorem \cite{Cercignani1988} gives the dissipation of entropy:
$$
\int_{\mathbb{R}^3} \mathcal{Q}(f)(v) \log (f)(v) d v \leq 0.
$$
Boltzmann’s H-theorem implies that any equilibrium distribution function has the form of a local Maxwellian distribution $\mathcal{M}_{\rho, u, T}$, i.e.,
\begin{equation}\label{eq:M}
\mathcal{Q}(f)=0 \quad \Leftrightarrow \quad f=\mathcal{M}_{\rho, u, T}:=\frac{\rho}{(2 \pi T)^{d_v / 2}} \exp \left(-\frac{|v-u|^2}{2 T}\right). 
\end{equation}
Here $\rho, u$ and $T$ are the density, bulk velocity and temperature of a gas, respectively:
\begin{equation}\label{eq:rho_u_T}
    \rho=\int_{v \in \mathbb{R}^{d_v}} f(v) d v, \  u=\frac{1}{\rho} \int_{v \in \mathbb{R}^{d_v}} v f(v) d v, \  T=\frac{1}{d_v \rho} \int_{v \in \mathbb{R}^{d_v}}|u-v|^2 f(v) d v.
\end{equation}

\subsection{Hydrodynamic limit}

When the Knudsen number $\varepsilon>0$ becomes very small, the macroscopic fluid dynamics, which describe the evolution of averaged quantities such as the density $\rho$, momentum $\rho u$ and temperature $T$ of the gas, by the compressible Euler or Navier-Stokes equations, become adequate.

Specifically, multiplying \eqref{eq:B_sto} by $m(v)$, integrating over $v$ and using the conservation property given in \eqref{eq:conservation}, one obtains the following local conservation law:
\begin{equation}\label{eq:moment}
\left\{\begin{array}{l}
\partial_t \rho+\nabla_x \cdot(\rho u)=0, \\
\partial_t(\rho u)+\nabla_x \cdot(\rho u \otimes u+\mathbb{P})=\mathbf{0}, \\
\partial_t E+\nabla_x \cdot(E u+\mathbb{P} u+\mathbb{Q})=0,
\end{array}\right.
\end{equation}
where the total energy $E$, the pressure tensor $\mathbb{P}$ and the heat flux $\mathbb{Q}$ are defined by
\begin{equation}\label{eq:energy}
\begin{array}{c}
        E=\displaystyle\int_{\mathbb{R}^{d_v}} \frac{1}{2}|v|^2 f d v=\frac{1}{2} \rho|u|^2+\frac{d_v}{2} \rho T,\\
\mathbb{P}=\displaystyle\int_{\mathbb{R}^{d_v}} (v-u) \otimes(v-u) f d v,\ 
\mathbb{Q}=\displaystyle\int_{\mathbb{R}^{d_v}} \frac{1}{2} (v-u)| v-\left.u\right|^2 f d v.
\end{array}\end{equation}
Note that the variables $\rho, u$ and $E$ in \eqref{eq:moment} depend on the random parameter $z$.

When $\varepsilon \to 0$, the distribution $f^{\varepsilon}$ formally converges to a Maxwellian $\mathcal{M}_{\rho, u, T}$ from \eqref{eq:M}, and thus the expression $\mathbb{P}$ and $\mathbb{Q}$ can be approximated by moments of the Maxwellian, i.e.,
$$
\mathbb{P}=p I, \quad \mathbb{Q}=0,
$$
where $p=\rho T$ is the pressure and $I$ is the identity matrix. Then the above local conservation law\eqref{eq:moment} reduces to the following compressible Euler system:
\begin{equation}\label{eq:E_sto}
\left\{\begin{array}{l}
\partial_t \rho+\nabla_x \cdot(\rho u)=0, \\
\partial_t(\rho u)+\nabla_x \cdot(\rho u \otimes u+\rho T I)=\mathbf{0}, \\
\partial_t E+\nabla_x \cdot((E+\rho T)u)=0 .
\end{array}\right.
\end{equation}
which is known as a zeroth order approximation with respect to $\varepsilon$ to the Boltzmann equation \eqref{eq:B_sto}.

By the Chapman-Enskog expansion, the compressible Navier-Stokes system give a first order approximation in $\varepsilon$ to the Boltzmann equation \cite{Cercignani1988, Filbet2015sisc}:
\begin{equation}\label{eq:NS}
\left\{\begin{array}{l}
\partial_t \rho+\nabla_x \cdot(\rho u)=0, \\
\partial_t(\rho u)+\nabla_x \cdot(\rho u \otimes u+\rho T I)=\varepsilon \nabla_x \cdot(\mu \mathbb{D}(u)), \\
\partial_t E+\nabla_x \cdot(u(E+\rho T))=\varepsilon \nabla_x \cdot\left(\mu \mathbb{D}(u) \cdot u+\kappa \nabla_x T\right),
\end{array}\right.
\end{equation}
where the traceless deformation tensor $\mathbb{D}(u)$ is given by 
$$
\mathbb{D}(u)=\nabla_x u+\left(\nabla_x u\right)^{\top}-\frac{2}{d_v}\left(\nabla_x \cdot u\right) I,
$$
the scalar quantities $\mu$ and $\kappa$ are viscosity and the thermal conductivity, respectively.

\section{An asymptotic preserving hybrid (APH) scheme}
\label{sec3:hybrid}

We briefly review the hybrid solver for the Boltzmann equation, as studied in \cite{Filbet2015sisc, Tiwari1998jcp}. The main idea involves a dynamic domain decomposition approach that partitions the computational domain into two distinct regions: kinetic cells governed by the Boltzmann equation and fluid cells described by macroscopic fluid equations. 
At each time step, the algorithm performs critical cell classification, evaluating all computational cells to update the respective kinetic and fluid index sets. 
In order to prevent wrong judgment (for example, treating a situation where the fluid is far from the thermal equilibrium as fluid cells), we implement dual criteria:  one identifying breakdowns of the hydrodynamic description (fluid$\rightarrow$kinetic transition), and another detecting when kinetic solutions approach local equilibrium (kinetic$\rightarrow$fluid transition). It is noted that computational efficiency is achieved through an optimized decomposition strategy that minimizes the kinetic subdomain size, thereby maximizing the utilization of the fluid solvers which need low computational cost. 

In this paper, we employ the AP scheme proposed by Filbet and Jin in \cite{Filbet2010jcp} for the kinetic cells, which features a BGK-type penalization of the nonlinear collision term. This approach maintains explicit solvability even with implicit temporal discretization and can capture the  macroscopic fluid dynamic limit for small Knudsen number. For the fluid cells, we utilize a conservative finite-volume TVD scheme  \cite{Toro2013springer}, chosen for its shock-capturing capabilities and numerical robustness. Developing particle-based acceleration techniques for the hybrid scheme will be one of our future investigation.
The forthcoming sections will provide detailed fluid$\rightarrow$kinetic and kinetic$\rightarrow$fluid transition criteria, followed by a complete description of the APH algorithm.

To enhance the regime transition accuracy, we need to simultaneously:
i) identify hydrodynamic description failures, ii) tracking distribution function's $L_v^1$-distance to local equilibrium at each time step. At a given time $t^n$, the space domain $D=D_f^n \cup D_k^n$ is decomposed in fluid grids $G_i \subset D_f^n$ described by the fluid field
$$
U_i^n:=\left(\rho_i^n, u_i^n, T_i^n\right) \simeq\left(\rho\left(t^n, x_i\right), u\left(t^n, x_i\right), T\left(t^n, x_i\right)\right), \ x_i \in G_i,
$$
and kinetic grids $G_j \subset D_k^n$ described by the particle distribution function
$$
f_j^n(v) \simeq f^\varepsilon\left(t^n, x_j, v\right), \  x_j \in K_j, \ v \in \mathbb{R}^{d_v}.
$$

Below we give the two criteria in detail, followed by the APH algorithm.

\textbf{Fluid$\rightarrow$Kinetic Criteria.} 
Considering the positive definite matrix \cite{Filbet2015sisc}
$$\mathcal{V}:=I+\bar{A}^{\varepsilon}-\frac{2}{3 \bar{C}^{\varepsilon}} \bar{B}^{\varepsilon} \otimes \bar{B}^{\varepsilon},$$
with the traceless matrix $\bar{A}^{\varepsilon}$ and the vector $\bar{B}^{\varepsilon}$ given by
$$
\begin{cases}\bar{A}^{\varepsilon}:=\dfrac{1}{\rho^{\varepsilon}} \int_{\mathbb{R}^d_v} A(V) f^{\varepsilon}(v) d v, & A(V)=V \otimes V-\dfrac{|V|^2}{d_v} I, \quad V(v)=\dfrac{v-u}{\sqrt{T}},\\ \bar{B}^{\varepsilon}:=\dfrac{1}{\rho^{\varepsilon}} \int_{\mathbb{R}^d_v} B(V) f^{\varepsilon}(v) d v, & B(V)=\dfrac{1}{2}\left[|V|^2-(d_v+2)\right] V, \quad V(v)=\dfrac{v-u}{\sqrt{T}},\end{cases}
$$
and $\bar{C}^{\varepsilon}$ the dimensionless fourth order moment of $f^\varepsilon$, given by
$$\bar{C}^{\varepsilon}:=\frac{2}{d_v \rho} \int_{\mathbb{R}^d_v}\left[\frac{|V|^2}{2}-\frac{d_v}{2}\right]^2 f^{\varepsilon}(v) d v,$$
the zeroth-order $\left(\varepsilon^0\right)$ compressible Euler system yields $\mathcal{V}_{\text {Euler }}=I$, while first-order $\left(\varepsilon^1\right)$ compressible Navier-Stokes system introduces dissipative terms:
\begin{equation}\label{matrix:V_NS}
    \mathcal{V}_{N S}:=\mathcal{V}_{\varepsilon}=I-\varepsilon \frac{\mu}{\rho T} D(u)-\varepsilon^2 \frac{2}{d_v} \frac{\kappa^2}{\rho^2 T^{d_v}} \nabla_x T \otimes \nabla_x T,
\end{equation}
where ($\rho, u, T$) are solution to the Navier-Stokes equations \eqref{eq:NS}.

Hence, the Euler approximation is valid only if the matrix $\mathcal{V}_{N S}$ is positive definite and  all eigenvalues $\lambda_{N S}$ of $\mathcal{V}_{N S}$ satisfy 
\begin{equation}\label{cri_FK}
    \left|\lambda_{N S}-1\right|>\eta_0, \quad \forall \lambda_{N S} \in \operatorname{Sp}\left(\mathcal{V}_{N S}\right),
\end{equation}
with $\eta_0$ a small parameter (we take $\eta_0=10^{-2}$ in simulations), ensuring $\mathcal{V}_{N S} \approx I$.

\textbf{Kinetic$\rightarrow$Fluid Criteria.} 
Here we use the $L_v^1$-distance between the distribution $f^{\varepsilon}$ and its Maxwellian $\mathcal{M}_{\rho, u, T}^{\varepsilon}$ to check if the system is locally at the thermodynamic equilibrium or not, thus the kinetic$\rightarrow$fluid transition criterion is given by
\begin{equation}\label{cri_KF1}
    \left\|f^{\varepsilon}(x, t, \cdot)-\mathcal{M}_{\rho, u, T}^{\varepsilon}(x, t, \cdot)\right\|_{L_v^1} \leq \delta_0,
\end{equation} 
with $\delta_0$ a small parameter (we take $\delta_0=10^{-4}$ in simulation). 
This aligns with the Chapman-Enskog expansion's small-remainder term $$\left\|\sum_{n>1} \varepsilon^i g^{(i)}(t,x,\cdot)\right\|_{L_v^1} \leq \delta_0$$ for the first-order hydrodynamic closure.

The complete APH algorithm is described in Algorithm \ref{alg:hybrid}. It is noted that at each time step, for the fluid grids processing, we first evaluate the fluid$\rightarrow$kinetic criterion and then update the fluid domain using the fluid solver; while for the kinetic grids processing, we initially advance the solution in the current kinetic domain via the kinetic solver and subsequently perform the criterion evaluation and update both domains.
This staggered approach guarantees that every grid is updated exactly once per iteration. To illustrate the necessity of this design, consider the alternative approach where domain updates are always performed prior to solver computations: such implementation could lead to pathological cases where certain grids oscillate indefinitely between kinetic and fluid domains without completing state updates. The algorithm avoids explicit loops over individual domain indices, instead employing vectorized operations to enhance computational efficiency. 

\begin{algorithm}[t]
\caption{An asymptotic preserving hybrid (APH) numerical algorithm}
\label{alg:hybrid}
\begin{algorithmic}[1]
\State \textbf{Input}: Initial domains $D_f^0$, $D_k^0$, fluid field $\{U_i^0\}$, particle distribution $\{f_j^0\}$
\State \textbf{Output}: Evolved solutions $\{U_i^{N+1}\}$, $\{f_j^{N+1}\}$, updated domains $D_f^{N+1}$, $D_k^{N+1}$

\For{$n=0$ to $N_{\text{steps}}$}  
    \Statex \textbf{ Fluid Grids ($G_i \subset D_f$) Processing:}
    \State Evaluate the fluid$\rightarrow$kinetic criterion \eqref{cri_FK} for each fluid grid $G_i\subset D_f$
    \State Define $D_{f\rightarrow k} := \{G_i \subset D_f | \text{ criterion \eqref{cri_FK} holds} \}$ 
    \State Update $D_f := D_f \setminus D_{f\rightarrow k}$ and $D_k := D_k \cup D_{f\rightarrow k}$
    \State Update the fluid field $U_i^{n+1}$ for the fluid domain $D_f$ using a fluid solver $U_i^{n+1}=\mathcal{F}_\Delta t(U_i^n),\  G_i\subset D_f$
    \State “Lift” the macroscopic field into the kinetic grid: $$f_i^n(v) := \mathcal{M}(\rho_i^n, u_i^n, T_i^n), \ G_i\subset D_{f\rightarrow k}$$  

    \Statex \textbf{ Kinetic Grids ($G_j \subset D_k$)  Processing:}
    \State Update the particle distribution $f_j^{n+1}$ for the kinetic domain $D_k$ using a kinetic solver $f_j^{n+1}=\mathcal{K}_{\Delta t}\left(f_j^n\right), \  G_j \subset D_k$
    \State Evaluate the kinetic$\rightarrow$fluid criterion \eqref{cri_KF1} for each kinetic grid $G_j\subset D_k$
    \State Define $D_{k\rightarrow f} := \{G_j \subset D_k | \text{ criterion \eqref{cri_KF1} holds} \}$ 
    
    \State Project the kinetic distribution to the macroscopic field: $$U_j^n := \int f_j^n\phi(v)dv\ \text{ for }\ \phi(v) \in \left(1,v,\frac{|v-u_j^n|^2}{d_v\rho_j^n}\right), \ G_j \subset D_{k\rightarrow f}$$
    
    \State Update $D_k := D_k \setminus D_{k\rightarrow f}$ and $D_f := D_f \cup D_{k\rightarrow f}$
\EndFor

\end{algorithmic}
\end{algorithm}

Finally, special attention must be paid to the consistent updating of ghost points when computing fluid fields for kinetic$\rightarrow$fluid transitions and particle distributions for fluid$\rightarrow$kinetic transitions, especially at boundary points and when adjacent grids are of different types.
In one dimensional spatial domain, considering the situation at time $t^n$ where the grids $G_{i-1}$ and $G_i$ are fluid, and the grids $G_{i+1}$ and $G_{i+2}$ are kinetic, we need to lift the hydrodynamic fields by setting
$$
f_{i-k}^n(v):=\mathcal{M}_{\rho_{i-k}^n, u_{i-k}^n, T_{i-k}^n}(v), \quad k=0,1, \quad \forall v \in \mathbb{R}^{d_v},
$$
to construct ghost distributions for kinetic evolution in $G_{i+1}$, and project the kinetic density by setting
$$
U_{i+k}^n:=\int_{\mathbb{R}^d} f_{i+k}^n\left(1, v, \frac{1}{d_v \rho_{i+k}^n}\left|v-u_{i+k}^n\right|^2\right) d v, \quad k=1,2.
$$
with flux matching to construct ghost moments for fluid evolution in $G_i$.

\section{Stochastic collocation based method}
\label{sec4:sC}

\subsection{Control variate multi-level Monte Carlo method}

By reviewing the standard Monte Carlo (MC), we present a control variate multilevel Monte Carlo (MLMC) method, where control variates are employed to quasi-optimize variance reduction locally using two subsequent levels \cite{Hu2021siamjuqMLMC,Hu2023jcp}.

\textbf{Standard Monte Carlo method. }Supposing that we generate $M$ independent and identically distributed (i.i.d.) random samples $f_0^i, i=1, \ldots, M$, according to the random initial condition $f_0(w ; x, v)$, then each sample $f_0^i(w ; x, v)$ gives rise to a unique solution $f^i(w ; x, v, t)$ to \eqref{eq:B_sto} at time $t$. From $f^i(w ; x, v, t)$, we can easily compute
$$ \rho^i(w ; x, t)=\int_{\mathbb{R}^3} f^i(w ; x, v, t) \mathrm{d} v, \quad u^i(w ; x, t)=\dfrac{1}{\rho^i}\int_{\mathbb{R}^3} v f^i(w ; x, v, t) \mathrm{d} v, $$
$$ T^i(w ; x, t)=\dfrac{1}{d_v \rho^i}\int_{\mathbb{R}^{d_v}} |v-u|^2 f^i(w ; x, v, t) \mathrm{d} v. $$
In the following, we use a single variable $q(w; x, t)$ to denote macroscopic quantities  $\rho(w; x, t)$, $u(w; x, t)$ or $T(w; x, t)$.

Given the samples $q_{\Delta x, \Delta t}^i, i=1, \ldots, M$, computed with mesh size $\Delta x$ and time step $\Delta t$ from initial data $f_0(w_i ; x, v)$ up to time $t$, the MC estimate of the expectation $\mathbb{E}[q(w ; x, t)]$ is now given by
\begin{equation}\label{eq:MC}
    \mathbb{E}[q(w ; x, t)] \approx E_M\left[q_{\Delta x, \Delta t}(w ; x, t)\right]:=\frac{1}{M} \sum_{i=1}^M q_{\Delta x, \Delta t}^i(w ; x, t) .
\end{equation}

\bigskip

\textbf{Quasi-optimal multilevel Monte Carlo method.}  
The MLMC method  is known as a multilevel discretization with each level $l$ owning $M_l$ ($l=0,1,\cdots,L-1$) samples in random space. In each level, we use different meshes in temporal and spatial discretizations for the deterministic solver. 
Suppose that we have $L$ levels of solutions $\left\{q_{\Delta x_l, \Delta t_l}\right\}_{l=0}^{L-1}$, from coarsest level $q_{\Delta x_0, \Delta t_0}$ to finest level $q_{\Delta x_{L-1}, \Delta t_{L-1}}$, then the control variate MLMC method is given by
\begin{equation}\label{eq:CVMLMC}
    \begin{aligned}
\mathbb{E}[q(w ; x, t)] & \approx E_{C V}^L\left[q_{\Delta x_{L-1}, \Delta t_{L-1}}(w ; x, t)\right] \\
& :=\prod_{i=0}^{L-1} \lambda_i E_{M_0}\left[q_{\Delta x_0, \Delta t_0}\right]+\sum_{l=1}^{L-1} \prod_{i=l}^{L-1} \lambda_i E_{M_l}\left[q_{\Delta x_l, \Delta t_l}-\lambda_{l-1} q_{\Delta x_{l-1}, \Delta t_{l-1}}\right],
\end{aligned}
\end{equation}
where $\lambda_l$ are the coefficients to be determined below. 
At each level $l$, $M_l$ i.i.d. samples $f_0^i, i=1, \ldots, M_l$ of the initial data $f_0$ are generated on meshes $\Delta x_l$ and $\Delta x_{l-1}$ respectively.

Denote $q_{(l)} = q_{\Delta x_l, \Delta t_l}, q^i_{(l)} = q_{\Delta x_l, \Delta t_l}^i, i=1,2,\ldots, M_l$, $l=0,1,\ldots,L-1$ and consider variance reduction for each pair of consecutive levels, known as the quasi-optimal MLMC method \cite{Hu2023jcp}, yields the estimate for $\left\{\lambda_l\right\}_{l=0}^{L-1}$:
$$
\begin{aligned}
&\lambda_{l-1}=\frac{\operatorname{Cov}\left[q_{(l)}, q_{(l-1)}\right]}{V\left[q_{(l-1)}\right]} \approx \frac{\sum_{i=1}^{M_l}\left(q_{(l)}^i-\bar{q}_{(l)}\right)\left(q_{(l-1)}^i-\bar{q}_{(l-1)}\right)}{\sum_{i=1}^{M_l}\left(q_{(l-1)}^i-\bar{q}_{(l-1)}\right)^2},\, l=0,1,\ldots,L-2,
\end{aligned}$$
and $\lambda_{L-1}=1$, where $\bar{q}_{(l)}=E_{M_l}\left[q_{(l)}\right]$. In simulations, the numerical solutions $q^i_{(l)}$ and $q^i_{(l-1)}$ are obtained by using the APH scheme for the Boltzmann equation \eqref{eq:B_sto} at time $t$.

\subsection{Multi-fidelity method}

In this section, we briefly review the multi-fidelity approximation, particularly bi- and tri-fidelity method studied in \cite{Narayan2014sisc,Zhu2014siamjuq}, followed by our choices of the high- and low-fidelity solvers in this project. 

\textbf{Multi-fidelity approximation. } 
The multi-fidelity algorithm for approximating the high-fidelity solution consists of offline and online stages. In the offline stage, we employ the cheap low-fidelity model to explore critical parameter points. Within the online stage, we construct multi-fidelity approximation by transferring low-fidelity approximation rules. The detailed algorithm is summarized in Algorithm \ref{alg:2}.
Specifically, the bi-fidelity approximation method employs two computational solvers: a more accurate but expensive high-fidelity solver and a faster but less precise low-fidelity solver. The tri-fidelity approach adds a medium-fidelity solver that balances cost and accuracy between these two extremes, providing more flexible trade-offs between computational efficiency and solution accuracy.
The two key implementation aspects in Algorithm \ref{alg:2} are point selection strategy in Step 3 and multi-fidelity approximation in Steps 5-6. For point selection strategy, see \cite{Narayan2014sisc,Zhu2014siamjuq} for more technical details and error analysis.

\textbf{The low- and high-fidelity solutions.} 
Let $u^H(z)$, $u^M(z)$ and $u^L(z)$ denote high-, medium- and low-fidelity solutions respectively, with $N$ the number of affordable low-fidelity simulation runs and $K \ll N$ the number of affordable high-fidelity simulation runs. 
For sample points $\gamma_K=\left\{z_{i_1}, \cdots, z_{i_K}\right\} \subset \Gamma_N$, define the low-fidelity approximation space: $$\mathscr{U}^{\mathcal{Y}}\left(\gamma_K\right)= \text{span} \{u^\mathcal{Y}(\gamma_K)\} = \text{span}
\{u^\mathcal{Y}(z_{i_1}), \cdots, u^\mathcal{Y}(z_{i_K})\},\  \mathcal{Y} = \{L,M\},$$
with $L$ and $M$ corresponding to the bi- and tri-fidelity method respectively,
and the high-fidelity approximation space:
$$\mathscr{U}^H(\gamma_K) = \text{span} \{u^H(\gamma_K)\} = \text{span}\{u^H(z_1),...,u^H(z_k)\}.$$

\begin{algorithm}[t]
	\caption{Bi-fidelity ($\mathcal{Y}=L, \mathcal{F}=B$) and tri-fidelity  ($\mathcal{Y}=M, \mathcal{F}=T$) approximation for a high-fidelity solution at given $z$.}	\label{alg:2}
	\begin{algorithmic}[1]
				\Statex \textbf{Offline:}
				\State Select a sample set $\Gamma_N=\left\{z_1, z_2, \ldots, z_N\right\} \subset I_z$.
				\State 
                Conduct simulation for the low-fidelity model $u^L(z_j)$ at each  $z_j \in \Gamma_N$. 
				\State 
                Select $K$ ``important" points from $\Gamma_N$ and construct the low-fidelity approximation space $\mathscr{U}^{\mathcal{Y}}\left(\gamma_K\right)= \text{span}
            \{u^\mathcal{Y}(z_{i_1}), \cdots, u^\mathcal{Y}(z_{i_K})\},\ \mathcal{Y} = \{L,M\}.$
				\State 
                Run the high-fidelity model at each sample point of $\gamma_K$ and
                construct the high-fidelity approximation space 
                $\mathscr{U}^H\left(\gamma_K\right)=\text{span} 
            \{u^H(z_{i_1}), \cdots, u^H(z_{i_K})\}.$
				\Statex \textbf{Online:}
				\State
                For any given $z \in I_Z$, compute the low-fidelity solution $u^{\mathcal{Y}}(z)$ and the corresponding projection coefficients: 
				\begin{align*}\label{Online-eq}
					u^{\mathcal{Y}}(z)\approx \mathcal{P}_{\mathscr{U}^{\mathcal{Y}}\left(\gamma_K\right)}
                \left[u^{\mathcal{Y}}(z)\right]= \sum_{k=1}^K c_k^{\mathcal{Y}}(z) u^{\mathcal{Y}}\left(z_{i_k}\right),  \ \mathcal{Y} = \{L, M\}.
			\end{align*}				
				\State 
                Construct the multi-fidelity approximation by applying the same approximation rule as in the low-fidelity model: 
				\begin{equation*}\label{uB-eqn}
					u^{\mathcal{F}}(z)=\sum_{k=1}^K c_k^{\mathcal{Y}}(z) u^H\left(z_{i_k}\right), \  \mathcal{F} = \{B, T\}, \  \mathcal{Y} = \{L, M\}.
		\end{equation*}
	\end{algorithmic}
\end{algorithm}

Rather than using the full distribution function $f$ as the high-fidelity solutions, we focus on its macroscopic quantities of interest:
$$u^H(\omega;\cdot)=\left[\rho^H(\omega;\cdot), u^H(\omega;\cdot), T^H(\omega;\cdot)\right]^{\top}$$
as the high-fidelity snapshot solutions. This choice ensures consistency with the macroscopic quantities computed from the low-fidelity model. The corresponding low-fidelity solutions computed from the Euler system by using \eqref{eq:E_sto} are: $$u^L(\omega;\cdot)=\left[\rho^L(\omega;\cdot), u^L(\omega;\cdot), T^L(\omega;\cdot)\right]^{\top}.$$ 
To ensure that the initial data for both fidelity models are consistent, the initial data of $\rho$, $u$ and $E$ are derived from the initial distribution $f_0$ for the Boltzmann equation through \eqref{eq:rho_u_T} and \eqref{eq:energy}.

\textbf{The low- and high-fidelity solvers.} 
To mitigate computational costs, we select the compressible Euler system \eqref{eq:E_sto} as our low-fidelity model. As a zeroth-order ($\varepsilon^0$) approximation to the Boltzmann equation in the fluid regime, the Euler system captures macroscopic variations with significantly lower computational effort both in runtime and memory because it has no dependence on velocity space.
So in this work, we shall employ a finite-volume scheme for the hydrodynamic limit system \eqref{eq:E_sto} as our low-fidelity solver and employ the APH scheme discussed in Section \ref{sec3:hybrid} for the Boltzmann equation \eqref{eq:B_sto}  as our high-fidelity solver.

\section{Numerical simulation}
\label{sec5:numerical}

In this section, we first conduct numerical validation of the hybrid method's accuracy and computational efficiency through two representative test cases spanning the kinetic-to-fluid regime transition in \Cref{sec:deter}. Building upon these deterministic verifications, we then introduce uncertainties by considering random perturbations in both initial conditions and collision operators in \Cref{sec:uncertain}. Specifically, we simulate two challenging scenarios: a discontinuous shock tube problem with uncertain initial conditions and a mixed-regime case with spatially varying Knudsen numbers with either uncertain initial conditions or uncertain collision parameters. For the last mixed-regime case, we implement and compare two advanced sampling techniques, i.e., the MLMC method and a multi-fidelity approach, to assess their respective capabilities in handling multiscale systems with uncertainties. 

We take $d_x=1$ in all simulations. Since the matrix \eqref{matrix:V_NS} for the Navier-Stokes case $(k=1)$ can be explicitly given by
$$
\mathcal{V}_{N S}=\left(\begin{array}{ccc}
1-\varepsilon \frac{\mu}{\rho T} \partial_x u^x-\varepsilon^2 \frac{\kappa^2}{\rho^2 T^3}\left(\partial_x T\right)^2 & 0 & 0 \\
0 & 1+\varepsilon \frac{\mu}{\rho T} \partial_x u^x & 0 \\
0 & 0 & 1+\varepsilon \frac{\mu}{\rho T} \partial_x u^x
\end{array}\right),
$$
we can then read its eigenvalues on its diagonal. The criterion \eqref{cri_FK} for the Euler case ($k=0$) for a fluid cell to be kinetic at the next iteration is then
\begin{equation}\label{cri_FK_2}
    \left|\varepsilon \frac{\mu}{\rho T} \partial_x u^x+\varepsilon^2 \frac{\kappa^2}{\rho^2 T^3}\left(\partial_x T\right)^2\right| > \eta_0 \quad \text { or } \quad\left|\varepsilon \frac{\mu}{\rho T} \partial_x u^x\right| > \eta_0.
\end{equation}
In our simulations, we employ \eqref{cri_FK_2} as the transition criterion: grids remain in fluid regime ($\mathcal{V}_{NS} \approx I$) only if both conditions fail; otherwise, they switch to kinetic grids. 

For the Boltzmann equation, an APH scheme is utilized \cite{Filbet2010jcp}. The collision term $\mathcal{Q}(f^n,f^n)$ in \eqref{eq:B_sto} is evaluated numerically using the fast spectral method developed in \cite{Mouhot2006mc}. 
For spatial discretization, we implement a second-order upwind MUSCL scheme with a minmod slope limiter \cite{Leveque1992springer}, which effectively suppresses possible oscillations near discontinuities or sharp gradients.

For the Euler system, a second-order total variation diminishing (TVD) scheme with a minmod slope limiter is applied \cite{Toro2013springer}.
Since the Euler system evolves macroscopic quantities rather than the distribution function $f$ of the Boltzmann equation, this approach achieves significantly reduced computational cost and memory requirements.

\subsection{Tests of the APH method for deterministic Boltzmann equation}
\label{sec:deter}

\subsubsection{Test One: Sod shock tube test}

This test examines the numerical solution of the $1 D_x \times 2 D_v$ Boltzmann equation \eqref{eq:B_sto} without uncertainties through the benchmark Sod shock tube problem. We compare the results of our method against a reference solution obtained by a high-resolution discretization of the full kinetic equation.

More precisely, the initial condition is given by the Maxwellian distribution 
$$
f^{i n}(x, v)=\mathcal{M}_{\rho(x), u(x), T(x)}(v), \quad \forall x \in [0,1], \quad v \in \mathbb{R}^2,
$$
computed from the following macroscopic quantities: 
$$
(\rho(x), u(x), T(x))= \begin{cases}(1,0,1) & \text { if } 0 \leq x \leq 0.5,\\ (0.125,0,0.25) & \text { if } 0.5<x \leq 1,\end{cases}
$$
under specular boundary conditions \cite{Cercignani1994}.

\begin{figure}[t]
	\centering
	\includegraphics[width=\textwidth]{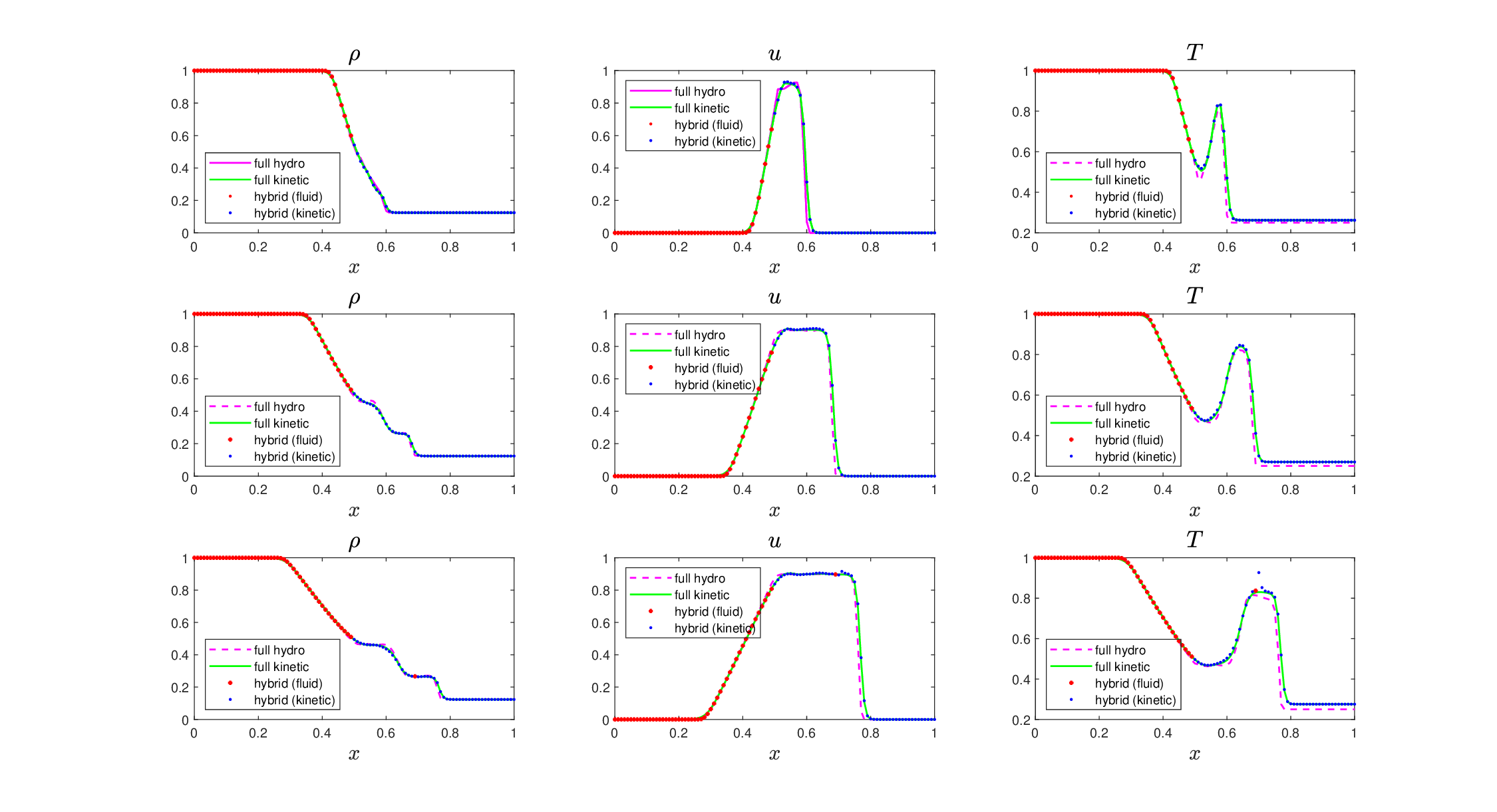}
	\caption{Test One: Sod shock tube test with $\varepsilon=10^{-4}$. Density $\rho$, mean
velocity $u$ and temperature $T$ at times $t = 0.05, 0.10 \text{ and } 0.15$.}\label{Ex1_Fig1}
\end{figure}

In Figure \ref{Ex1_Fig1}, we present the numerical results for the Sod shock tube test using the APH method with $\varepsilon=10^{-4}$. For the kinetic regime, we employ a fast spectral method for velocity discretization of the collision operator \cite{Pareschi2000sinum}, with $N_x=100$ in the spatial domain and $N_v=32$ in each velocity dimension within a truncated domain  $[-8,8]^2$. Temporal discretization is set with $\Delta t=8\times 10^{-4}$. For the Euler regime, we employ a classical TVD scheme  \cite{Toro2013springer} with the same spatial and temporal discretization for the kinetic regime. We present several snapshots of the density, mean velocity and temperature at different times $t=0.1$, $t=0.10$ and $t=0.15$. The solution by the full kinetic scheme, the full fluid scheme and the APH scheme are shown in Figure \ref{Ex1_Fig1}. Notably, the APH method achieves excellent agreement with the reference solution (the full kinetic scheme with $N_x=400$ mesh sizes) for density, mean velocity, and temperature, with only minor local discrepancies. 

The computational efficiency of the hybrid scheme is particularly evident in this setting. For the same configuration ($N_x=100, N_v=32$), the APH method runs 1.97 times faster than the full kinetic model. Here, the APH scheme maintains high accuracy while significantly reducing computational cost, demonstrating its suitability for high-dimensional simulations where efficiency gains will be even more substantial.

\subsubsection{Test Two: Blast wave}

We now consider the blast wave test case using initial conditions defined by a Maxwellian distribution 
$$
f^{i n}(x, v)=\mathcal{M}_{\rho(x), u(x), T(x)}(v), \quad \forall x \in[-0.5,0.5], \quad v \in \mathbb{R}^2,
$$
with
$$
(\rho(x), u(x), T(x))= \begin{cases}(1,1,2) & \text { if } x<-0.3, \\ (1,0,0.25) & \text { if }-0.3 \leq x \leq 0.3, \\ (1,-1,2) & \text { if } x \geq 0.3,\end{cases}
$$
over the spatial domain $x \in[-0.5,0.5]$ and $v \in \mathbb{R}^3$, under periodic boundary conditions. 

In this test, we perform simulations across a range of Knudsen numbers spanning from the rarefied regime to the fluid limit, presenting numerical results for $\varepsilon=10^{-1}$ and $10^{-4}$ cases. For both cases, the APH scheme employs $N_x=100$ spatial points and a $32 \times 32$ velocity grid over a truncated velocity domain $[-8,8]^2$, with results compared against full kinetic and hydrodynamic limit solutions.

\begin{figure}[t]
	\centering
	\includegraphics[width=\textwidth]{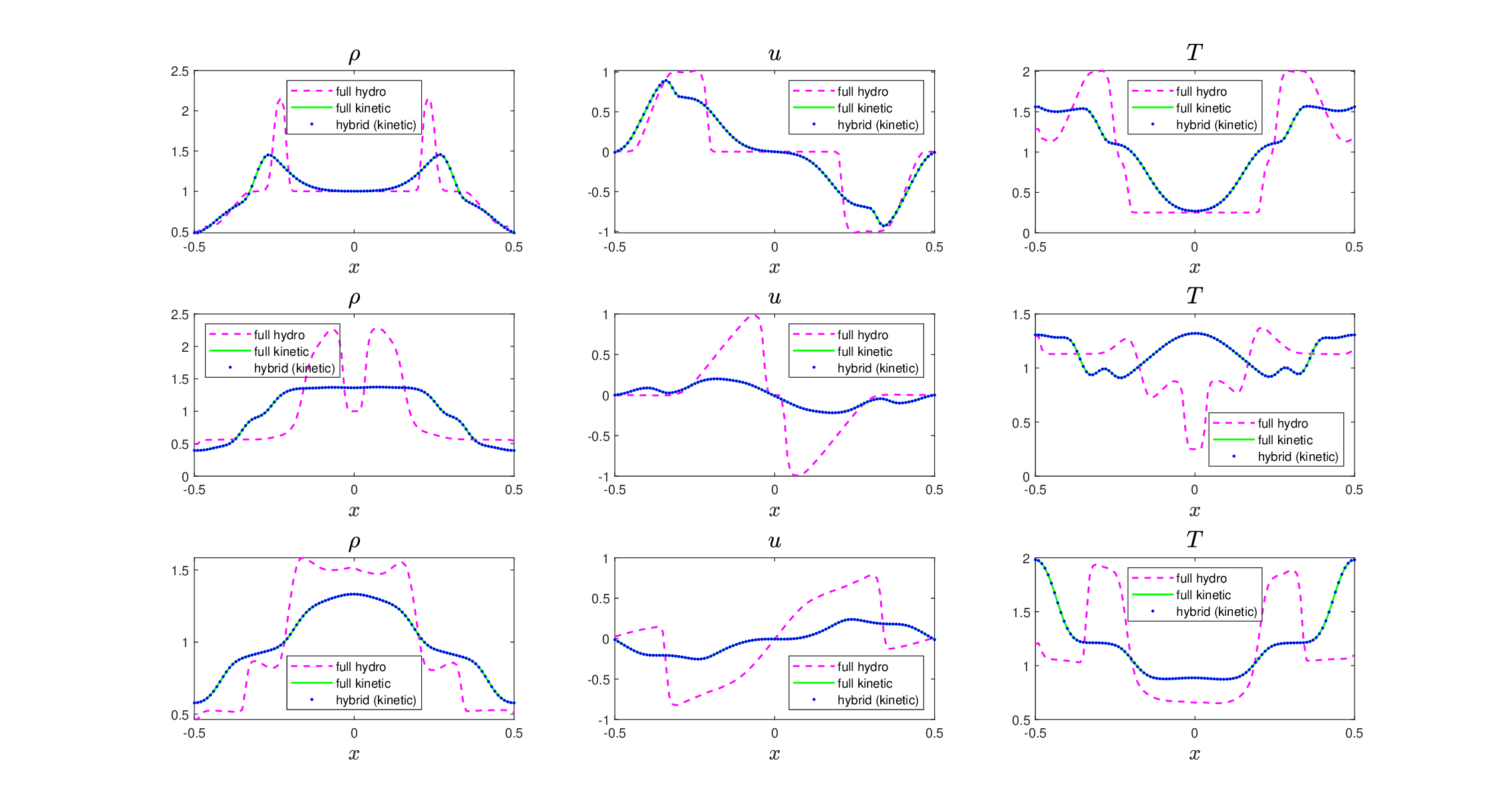}
	\caption{Test Two: Blast wave with  $\varepsilon=0.1$. Density $\rho$, mean
velocity $u$ and temperature $T$ at times $t = 0.05, 0.10 \text{ and } 0.35$.}\label{Ex2_Fig1eps1}
\end{figure}

\begin{figure}[t]
	\centering
	\includegraphics[width=\textwidth]{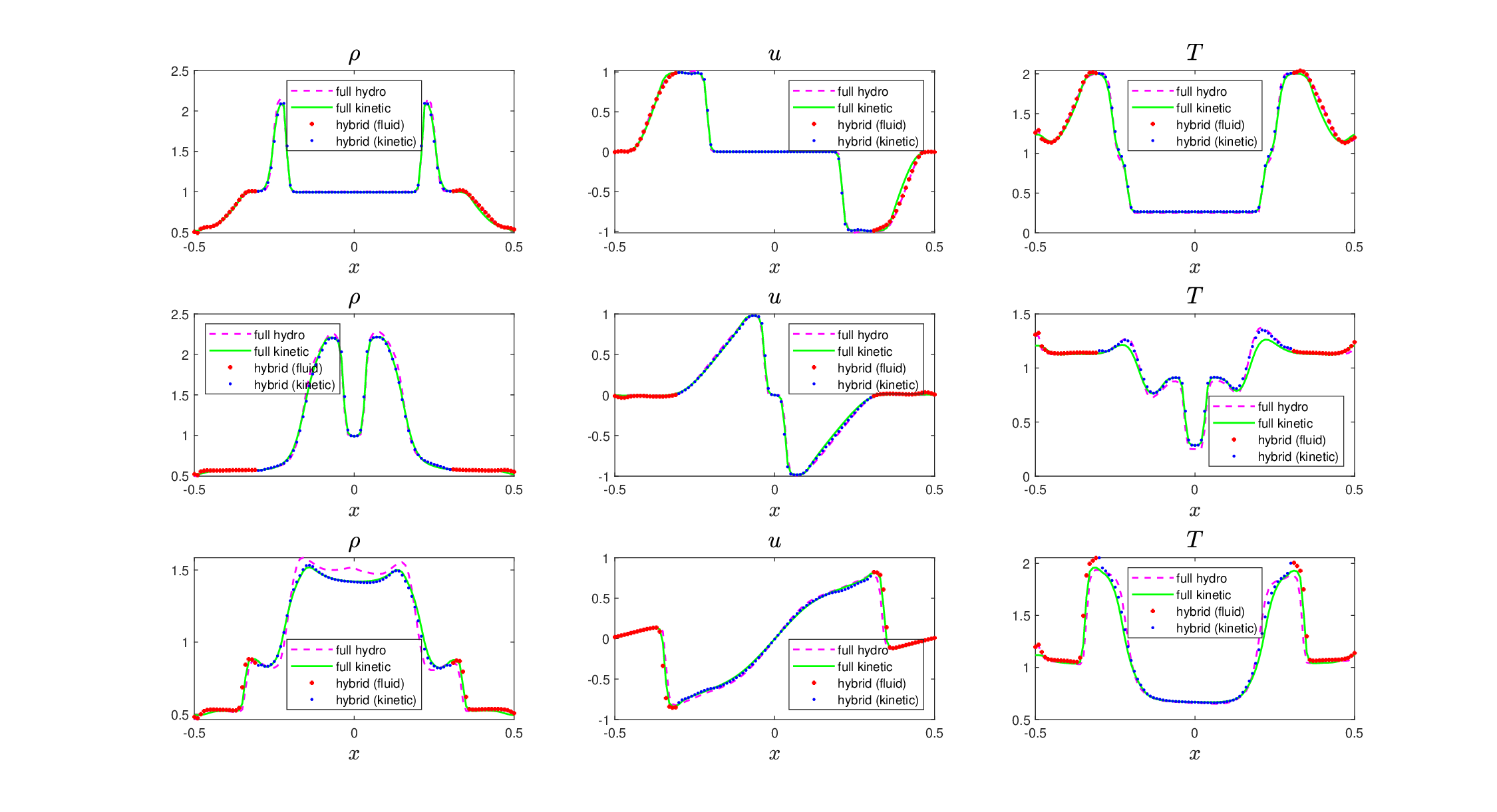}
	\caption{Test Two: Blast wave with  $\varepsilon=10^{-4}$. Density $\rho$, mean
velocity $u$ and temperature $T$ at times $t = 0.05, 0.10 \text{ and }  0.35$.}\label{Ex2_Fig2eps4}
\end{figure}

For the rarefied regime $(\varepsilon=0.1)$ shown in Figure
\ref{Ex2_Fig1eps1}, we observe that while the fluid solution is far from the kinetic one as expected for large Knudsen numbers, the hybrid method demonstrates excellent performance by accurately identifying and resolving non-equilibrium zones, yielding results closely aligned with the kinetic solution even at large times.

In the near-fluid regime ($\varepsilon=10^{-4}$) depicted in Figure \ref{Ex2_Fig2eps4}, the zeroth-order hybrid scheme achieves a good agreement with the reference solution while confining kinetic computations to minimal regions.  We observe that although the purely macroscopic model fails to provide accurate results for large time, the hybrid method consistently matches the full kinetic solution, demonstrating a $1.51$ computational speedup over the full kinetic solver. This efficiency gain becomes particularly significant when considering extensions to higher-dimensional spatial and velocity spaces, where kinetic computations are prohibitively expensive. 
When applying the first-order correction (Compressible Navier-Stokes solver), the error for the fluid solver will be further reduced (see \cite{Filbet2015sisc}).

In conclusion, the APH method dynamically adjusts to a range of Knudsen numbers, enabling accurate resolution in both rarefied and near-fluid regimes. Computational savings would scale further in high-dimensional settings and benefit for kinetic problems with uncertainties. Therefore maintaining accuracy and computational efficiency underscores the potential of the APH method for broader applications in multiscale flow simulations.

\subsection{Tests for Boltzmann equation with uncertainties}
\label{sec:uncertain}

In this section, we present numerical results for the Boltzmann equation \eqref{eq:B_sto} with either random initial conditions or random collision operators. 

First, extending our analysis of the Sod shock tube benchmark (Test One), we incorporate uncertainty quantification for the Boltzmann equation through the MLMC method. 
In \Cref{sec:Test3}, we consider the uncertainty in the initial conditions at a fixed Knudsen number $\varepsilon=10^{-4}$ and mainly focus on the impact of spatial resolution for such cases with discontinuous initial conditions.
\Cref{sec:Test4} investigates uncertainty propagation in the Boltzmann equation under two types of uncertainties, i.e., perturbations in the initial conditions and variability in the collision operator parameters,  which allows us to assess how different types of uncertainties propagate through the multiscale system and influence the numerical solutions.
We also conduct the mixed regime case with random initial conditions by two advanced UQ techniques: the MLMC method and the multi-fidelity method. By quantifying and comparing the computational cost versus accuracy for these methods, we finally provide guidelines for their optimal application in kinetic theory problems.

As the solution constitutes a random field, the numerical error is a random quantity as well. To quantify the approximation error, we conduct $J=10$ independent simulations for each method, obtaining the corresponding approximations $\left\{q^{(j)}(x, t)\right\}_{j=1}^J$ for the macroscopic quantities $q \in\{\rho, u, T\}$ and compute their statistical averages.

Two error quantification approaches are evaluated: one is the global error defined in the $L^2\left(\Omega ; L^2(D)\right)$ norm by:
\begin{equation}\label{err:global}
\text{Err}_q(t)=\left(\frac{1}{J} \sum_{j=1}^J\left\|q^{(j)}(\cdot, t)-q_{\mathrm{ref}}(\cdot, t)\right\|_{L^2(D)}^2\right)^{1/2},
\end{equation}
where $\|\cdot\|_{L^2(D)}$ denotes the $L^2$-norm over the physical domain, 
and the other is the pointwise error defined in spatial domain by:
\begin{equation}\label{err:pointwise}
\text{Err}_q(x, t)=\left(\frac{1}{J} \sum_{j=1}^J\left(q^{(j)}(x, t)-q_{\mathrm{ref}}(x, t)\right)^2\right)^{1/2}.
\end{equation}
Here, $q_{\text {ref }}$ denotes the reference solution computed using a high-resolution kinetic solver with $N_x=400$ spatial grid points.

\subsubsection{Test three:  Test One with uncertainty in initial data} 
\label{sec:Test3}

The Sod shock tube test with uncertainties presents a challenging problem due to its discontinuous initial conditions, which complicate both uncertainty propagation and numerical resolution. Assume the random initial distribution
$$
f_0=\dfrac{\rho_0}{2 \pi T_0} e^{-\frac{\left|v-u_0\right|^2}{2 T_0}} \quad \forall x \in[0,1], \quad v \in \mathbb{R}^2,
$$
where the initial data $\rho_0$, $u_0$ and $T_0$ are given by
$$
\left\{\begin{array}{llll}
\rho=1, & u=0, & T=1+0.4 \sum_{k=1}^{d_1} \frac{z_k^T}{2 k}, & x \leq 0.5, \\
\rho=\frac{1}{8}, & u=0, & T=\frac{1}{4}\left(1+0.4 \sum_{k=1}^{d_1} \frac{z_k^T}{2 k}\right), & x>0.5 .
\end{array}\right.
$$
Here the random variable $\mathbf{z}^T=\left(z_1^T, \cdots, z_{d_1}^T\right)$ in the initial temperature obeys the uniform distribution on $[-1,1]$. Set the random space dimension $d_1=5$
and the Knudsen number $\varepsilon=10^{-4}$.
The computational domain is set as $x \in[0,1]$ in physical space and $v \in[-8,8]^2$ in velocity space. The simulation employs the same discretization parameters as Test One. All the results are obtained using the quasi-optimal MLMC method \eqref{eq:CVMLMC}.

\begin{figure}[t]
	\centering
	\includegraphics[width=0.9\textwidth]{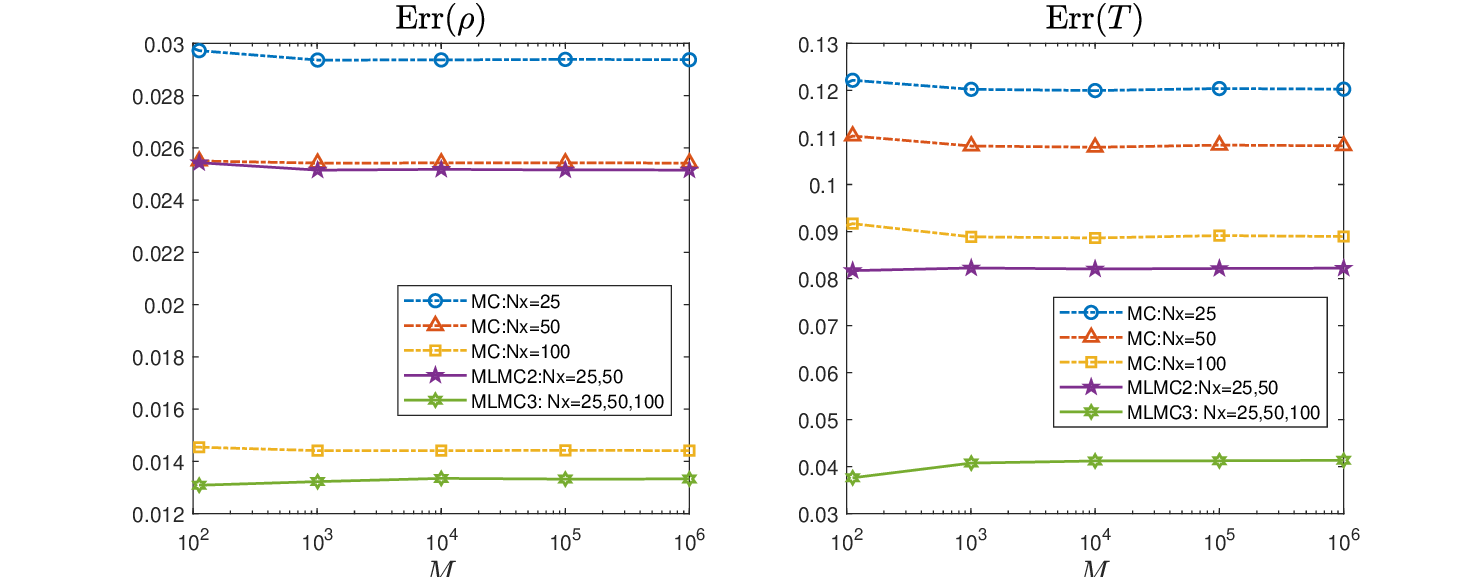}
	\caption{Test Three: shock wave with uncertainty in initial condition. Error of Density $\rho$ and temperature $T$ at times t = 0.05.}\label{Ex1_errI_MLMC_EPS4}
\end{figure}

To demonstrate the accuracy and computational efficiency of the MLMC method, we perform comparative studies using the following spatial discretization schemes:
For the standard Monte Carlo (MC) method, we employ two resolution levels with $N=50$ and $N=100$ grid points respectively. These are compared against two MLMC implementations:
\begin{itemize}
    \item a two-level MLMC scheme with mesh sizes $N_0=25$ and $N_1=50$;
    \item a three-level MLMC scheme with hierarchically refined meshes $N_0=25$, $N_1=50$ and $N_2=100$.
\end{itemize}
Figure \ref{Ex1_errI_MLMC_EPS4} presents the global errors \eqref{err:global} for macroscopic quantities $\rho$ and $T$ at $t=0.005$ under varying sample sizes $M_0=112$ to $10^6$. Due to the solution's low regularity, the errors exhibit early saturation - while statistical errors dominate for small sample sizes, spatial discretization errors become predominant as sample sizes increase. The results demonstrate clear advantages of the MLMC approach: 
\begin{itemize}
    \item [(1)] The two-level MLMC implementation (with mesh sizes $N_0=25, N_1=50$) demonstrates better accuracy than standard MC using $N=50$ grid points at the same level of sample size; 
    \item [(2)] The three-level MLMC ($N_0=25, N_1=50, N_2=100$) using only $M_0/16$ samples at the finest level outperforms both MC with $N=100$ using $M_0$ samples and the two-level MLMC ($N_0=25, N_1=50$) implementation.
\end{itemize}
 Therefore, MLMC maintains the performance advantage despite the reduced sample counts at finer levels ($M_1=M_0/4, M_2=M_0/16$), demonstrating how its hierarchical sampling strategy achieves computational efficiency while maintaining accuracy.

\begin{figure}[t]
	\centering
	\includegraphics[width=\textwidth]{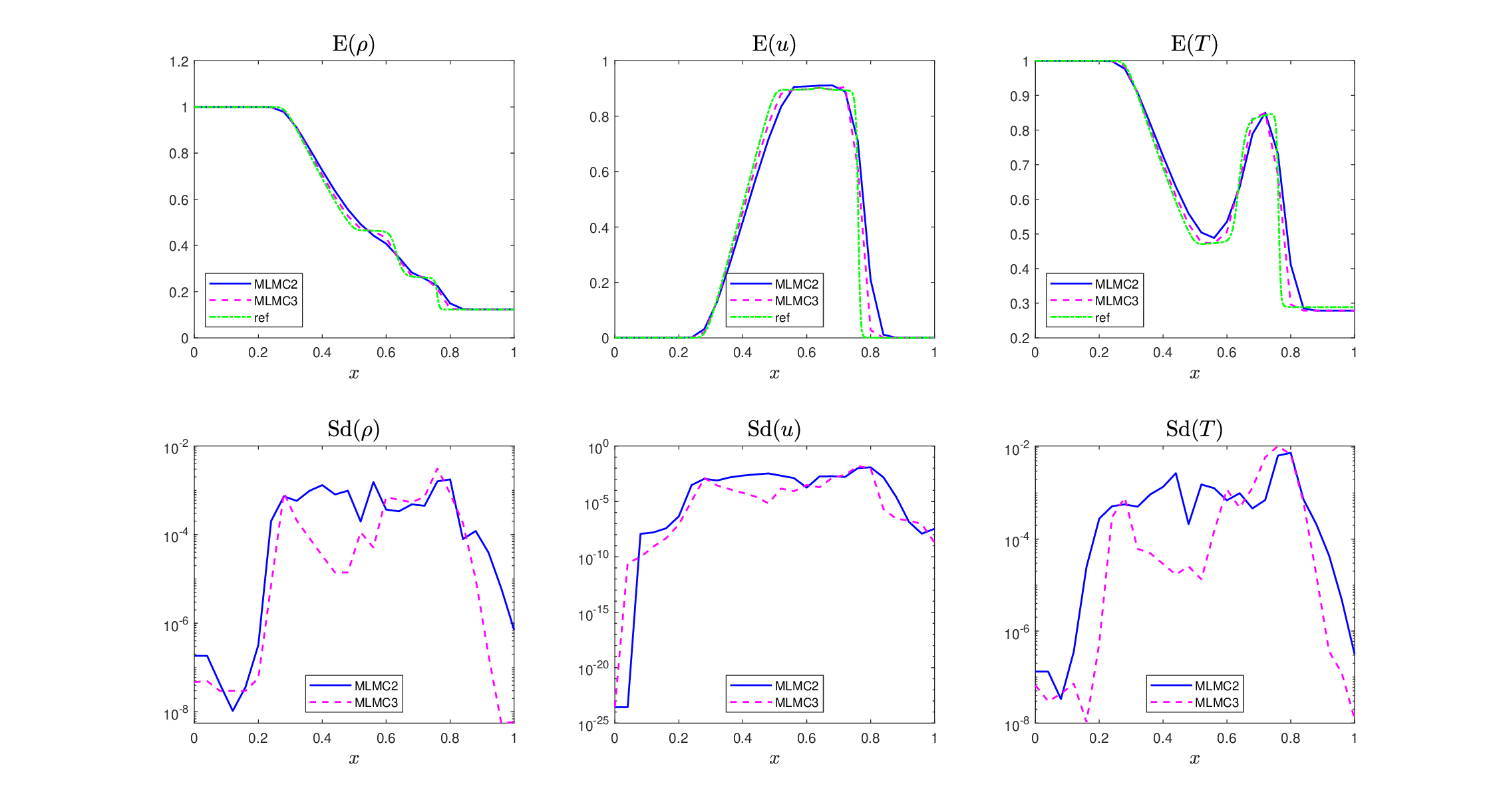}
	\includegraphics[width=\textwidth]{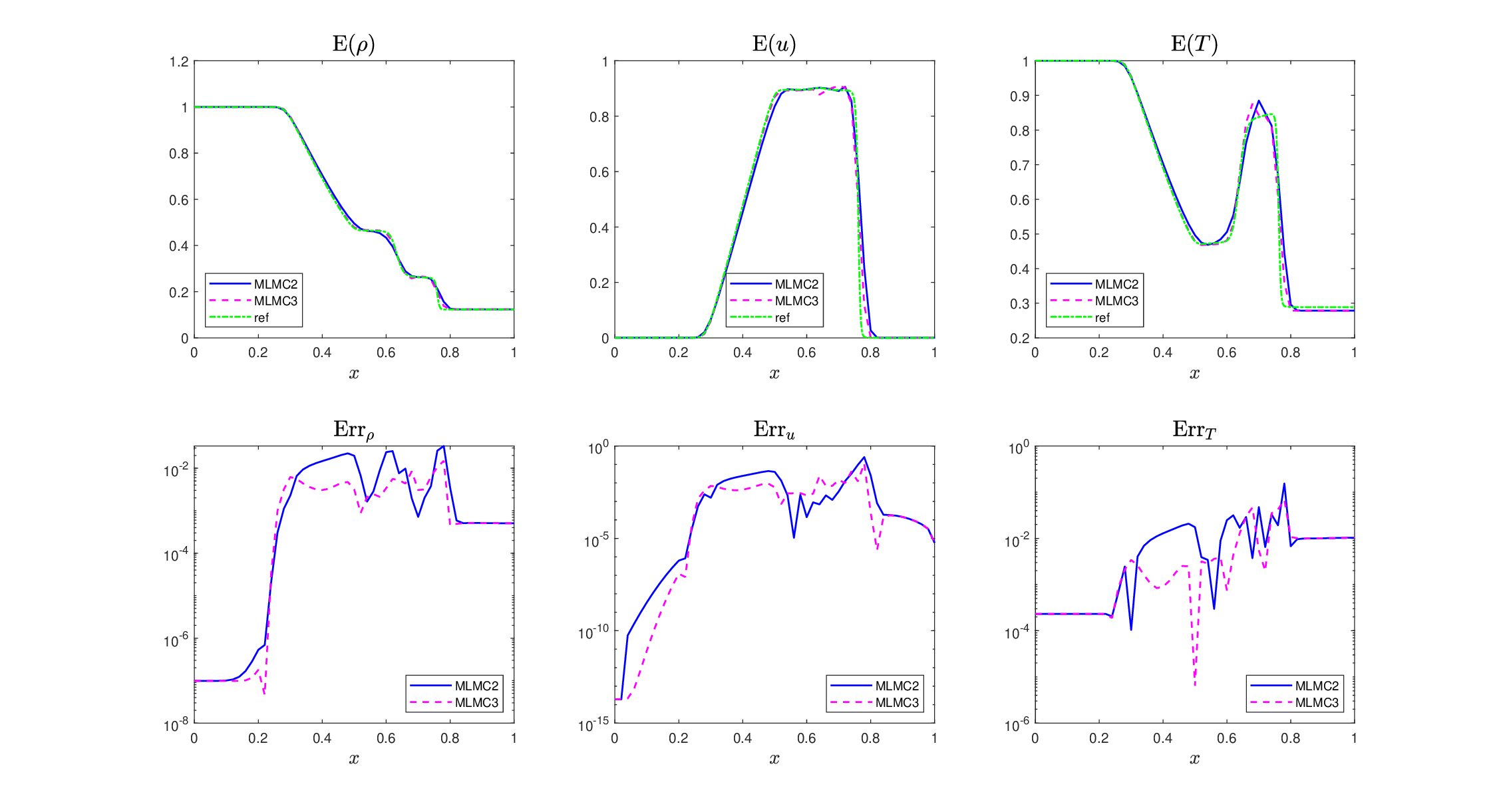}
	\caption{Test Three: shock wave with uncertainty in initial condition. The first two rows: Expectation and standard deviation of density $\rho$, mean velocity $u$ and temperature $T$ at time t = 0.15 by two-level MLMC ($N_0 = 25, N_1=50$) and three-level MLMC ($N_0 = 25, N_1=50, N_2=100)$.
    The last two rows: Expectation and pointwise error of density $\rho$, mean velocity $u$ and temperature $T$ at time t = 0.15 by two-level MLMC ($N_0 = 50, N_1=100$) and three-level MLMC ($N_0 = 50, N_1=100, N_2=200)$.}\label{Ex1_I015_MLMC100_EPS4}
\end{figure}

The impact of spatial resolution is examined in Figures \ref{Ex1_I015_MLMC100_EPS4} at $t=0.15$. With coarse grids ($N_0=25, N_1=50$ for two-level; $N_0=25$, $N_1=25, N_1=100$ for three-level), the first two rows in Figure \ref{Ex1_I015_MLMC100_EPS4} shows insufficient solution capture, though the three-level MLMC still yields smaller variance than its two-level counterpart. The variance is approximated via its statistical formulation, calculating the two expectations $\mathbb{E}\left[q^2\right]$ and $\mathbb{E}[q]$ separately and using them to obtain $\mathbb{V}[q]=\mathbb{E}\left[q^2\right]-(\mathbb{E}[q])^2$. 
 When refined to $N_0=50, N_1=100$ (two-level) and $N_0=50, N_1=100, N_2=200$ (three-level), the last two rows in Figure \ref{Ex1_I015_MLMC100_EPS4} demonstrates significantly improved accuracy in all macroscopic quantities $\rho, u, T$, quantified through our defined pointwise error \eqref{err:pointwise}. The results clearly show that the three-level MLMC method outperform the two-level MLMC method, and the MLMC method with finer meshes outperform the MLMC method with coarser meshes, especially in regions where the solution presents strong variations, namely close to the shock position.

On the one hand, the results demonstrate MLMC's key advantages over standard MC: (1) MLMC's hierarchical structure provides better error control than MC at equivalent computational cost, (2) three-level implementations achieve greater accuracy and variance reduction than two-level versions, and (3) grid refinement is crucial for resolving low-regularity solutions, especially near discontinuities.  These findings collectively establish MLMC's superiority over standard MC for this class of problems, particularly when handling discontinuous problems with limited computational resources.

On the other hand, the results shows that the APH-based MLMC method demonstrates robust performance across regimes, with particular advantages for lower-regularity solutions where finer meshes are needed. While theoretically converging to the exact solution with sufficient mesh refinement, in practice the control variate MLMC shows particularly strong performance in capturing solution features even with moderate mesh sizes.

\subsubsection{Test Four: mixed regime}
\label{sec:Test4}

\begin{figure}[t]
	\centering
	\includegraphics[width=0.4\textwidth]{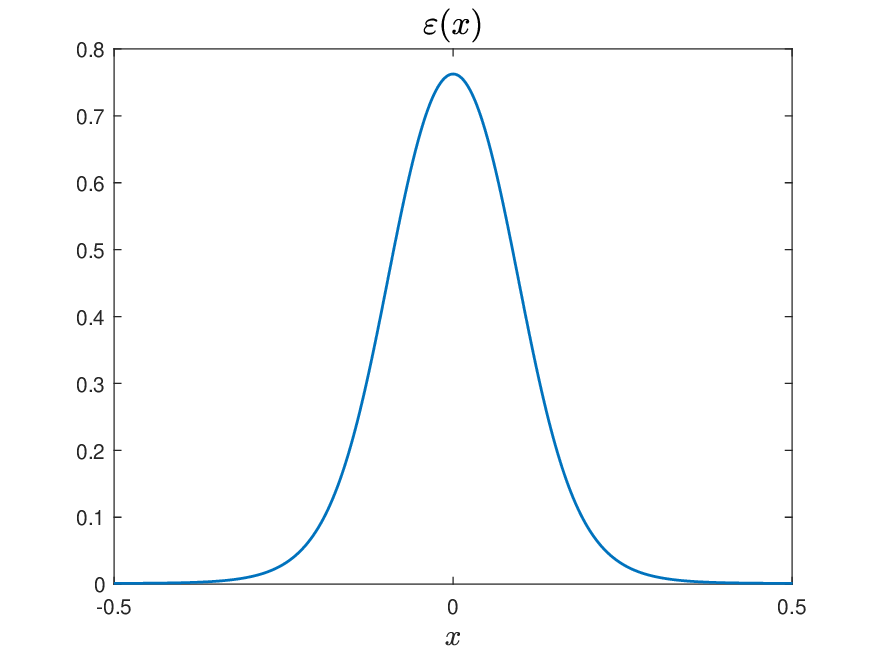}
	\caption{$\varepsilon(x)$.}\label{Fig:epsilon_x}
\end{figure}

 Now we consider the Boltzmann equation \eqref{eq:B_sto} with the Knudsen number $\varepsilon>0$ depending on the space variable in a wide range of mixing scales.
In this problem, $\varepsilon: \mathbb{R} \mapsto \mathbb{R}^{+}$is given by
$$
\varepsilon(x)=\varepsilon_0 +\frac{1}{2}[\tanh (1-11 x)+\tanh (1+11 x)],
$$
which varies smoothly from $\varepsilon_0$ to $O(1)$. Set $\varepsilon=10^{-3}$ in simulations, as shown in Figure \ref{Fig:epsilon_x}.

This numerical test is difficult because different scales are involved. It requires a good accuracy of the numerical scheme for all range of $\varepsilon$. In order to focus on the multi-scale nature, here we consider periodic boundary conditions. Furthermore, to increase the difficulty, let the initial distribution $f_0$ follow a double-peak non-equilibrium initial data \cite{Filbet2010jcp}, i.e., consider an initial distribution which is far from the local equilibrium of the collision operator:
$$
f_0=\frac{\rho_0}{2}\left[\exp \left(-\frac{\left|v-u_0\right|^2}{T_0}\right)+\exp \left(-\frac{\left|v+u_0\right|^2}{T_0}\right)\right], \quad x \in[-0.5, 0.5],\  v \in \mathbb{R}^2,
$$
with $u_0=(3 / 4,-3 / 4)$,
$$
\rho_0=\frac{2+\sin (2 \pi x)}{2}\left(1+0.4 \sum_{k=1}^{d_1} \frac{z_k^\rho}{2 k}\right), \quad T_0=\frac{5+2 \cos (2 \pi x)}{20}\left(1+0.4 \sum_{k=1}^{d_1} \frac{z_k^T}{2 k}\right).
$$
Here $\mathbf{z}^\rho=\left(z_1^\rho, \cdots, z_{d_1}^\rho\right), \mathbf{z}^T=\left(z_1^T, \cdots, z_{d_1}^T\right)$ represent the random variables in the initial density and temperature. In this test, we simulate three cases: 
\begin{itemize}
\item [\textbf{(a)}] Mixed regimes with uncertainty in initial conditions using the APH-based MLMC method;
\item [\textbf{(b)}] Mixed regimes with uncertainty in collision operators using the  APH-based MLMC method;
\item [\textbf{(c)}] Mixed regimes with uncertainty in initial conditions using the  APH-based multi-fidelity method.
\end{itemize}
In Test Four (a) and Test Four (c), set $d_1=7$, thus this is a $d=14$ dimensional problem in the random space.
In Test Four (b), uncertainty in collision operators is considered.

\begin{figure}[t]
	\centering
	\includegraphics[width=0.9\textwidth]{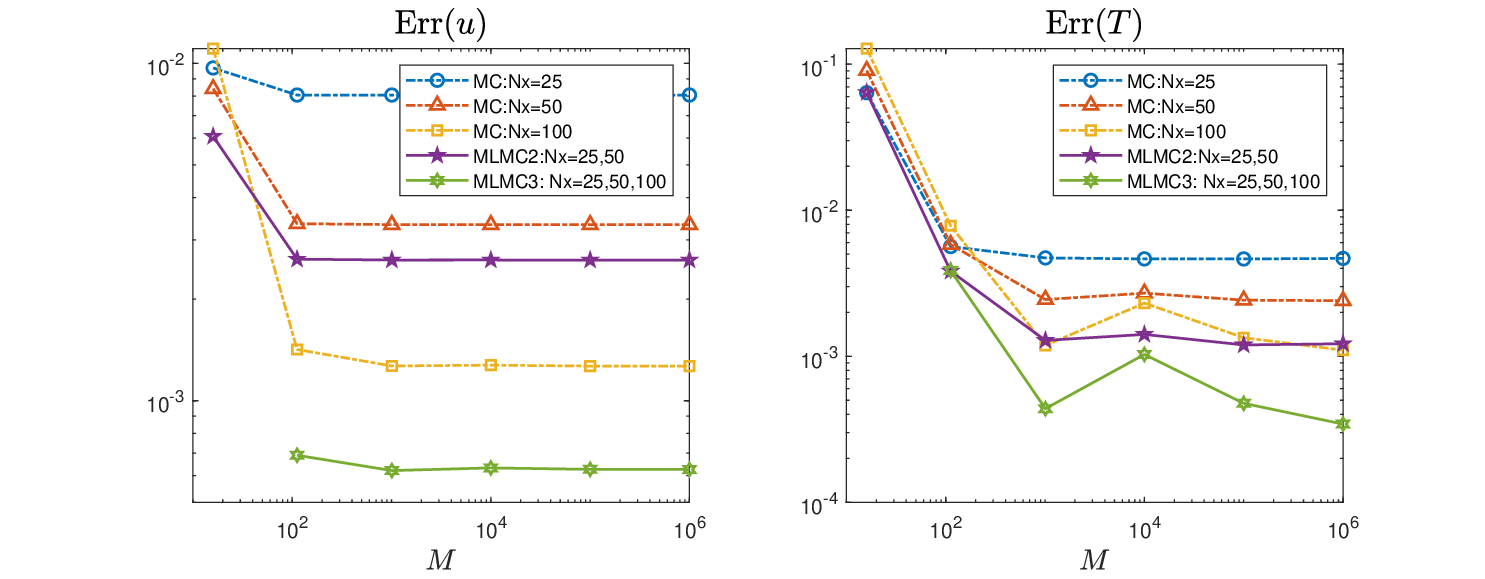}
	\caption{Test Four. Error of mean velocity $u$ and temperature $T$ at time $t=0.02$ by MC and MLMC methods.}\label{Fig:Ex5_err_MC_MLMC}
\end{figure}

\textbf{Test Four (a).} 
Similar as the previous example, we plot the errors \eqref{err:global} by MC and MLMC methods in Figure \ref{Fig:Ex5_err_MC_MLMC} as those in Figure \ref{Ex1_errI_MLMC_EPS4}. Due to better solution smoothness in this test case, the MLMC method achieves remarkable accuracy (up to at least $10^{-3}$) even with coarse spatial discretization ($N_0=25$). The results demonstrate that: 
\begin{itemize}
    \item [(1)] three-level MLMC ($N_0=25, N_1=50, N_2=100$) provides higher accuracy than two-level MLMC ($N_0=25, N_1=50$) at equivalent computational cost;
    \item [(2)] MLMC consistently outperforms the standard MC with the same level of sample sizes.
\end{itemize}

\begin{figure}[t]
	\centering
	\includegraphics[width=0.9\textwidth] {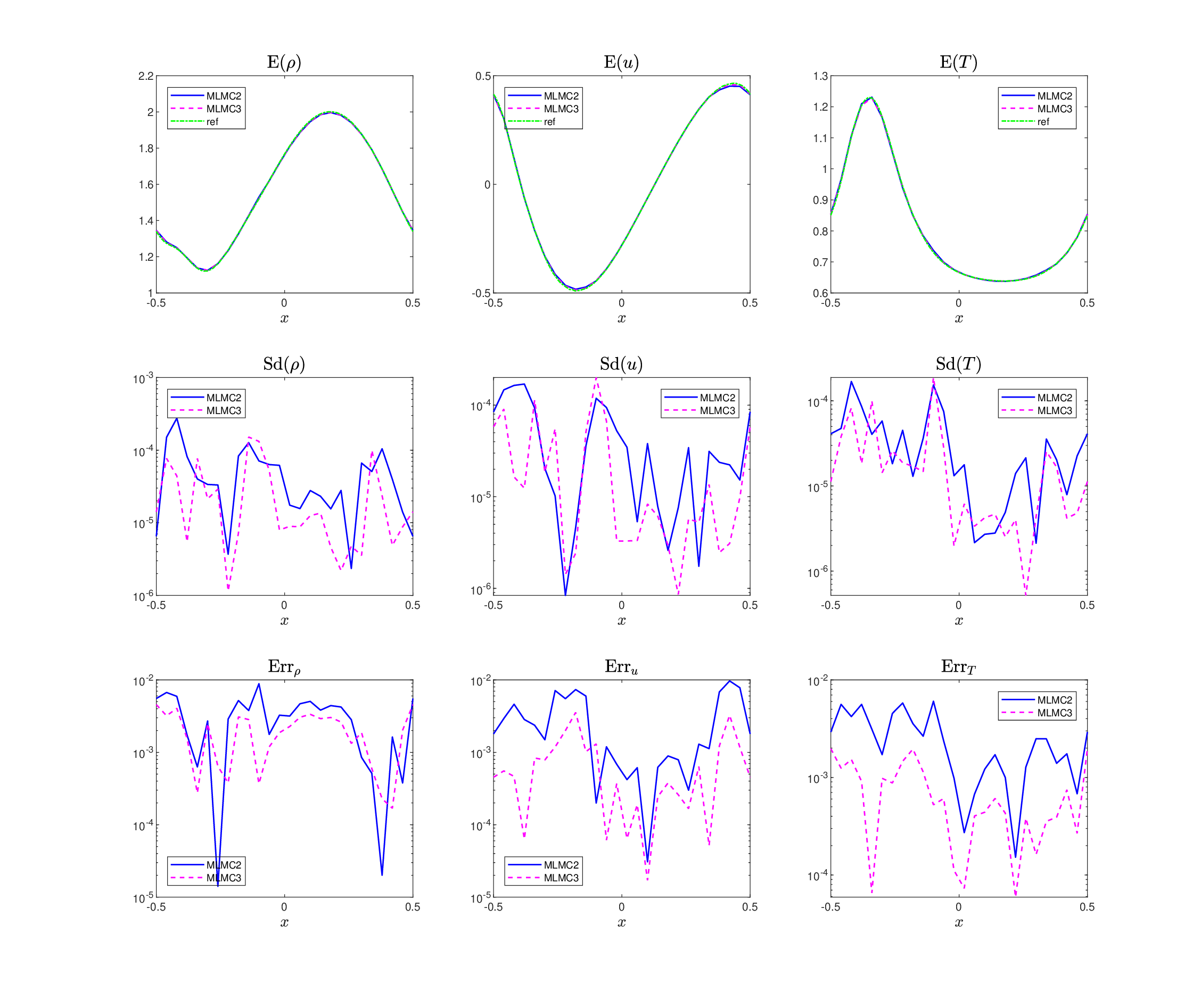}
	\caption{Test Four (a): Mixed regime with uncertainty in initial condition. Expectation, standard deviation and pointwise error of density $\rho$, mean velocity $u$ and temperature $T$ at time $t=0.2$ by the two-level MLMC ($N_0=25, N_1=50$) and the three-level MLMC ($N_0=25, N_1=50, N_2=100$).}\label{Fig:Ex5_MLMC3_100}
\end{figure}

In Figure \ref{Fig:Ex5_MLMC3_100}, we present the expectation, standard deviation and pointwise error of density $\rho$, mean velocity $u$ and temperature $T$ at time $t=0.2$ by the two-level MLMC ($N_0=25, N_1=50$) and the three-level MLMC ($N_0=25, N_1=50, N_2=100$). All reference solutions are computed via the full kinetic solver with $N_x=400$ spatial meshes. We observe that the APH-based MLMC method performs effectively with only 25 zeroth-level grids, benefiting from solution smoothness despite spatially varying $\varepsilon$.

This test validates the robustness of the APH-based MLMC method for mixed-regime problems with uncertain initial conditions, particularly demonstrating its ability to maintain accuracy while reducing computational costs through grid hierarchies.

\textbf{Test Four (b).} Assume the uncertain collision kernel in the form
$$
b\left(z^b\right)=1+0.5 z_1^b.
$$
Here $\mathbf{z}^b=z_1^b$ represent the random variables in the colllision kernel.

\begin{figure}[t]
	\centering
	\includegraphics[width=0.9\textwidth]{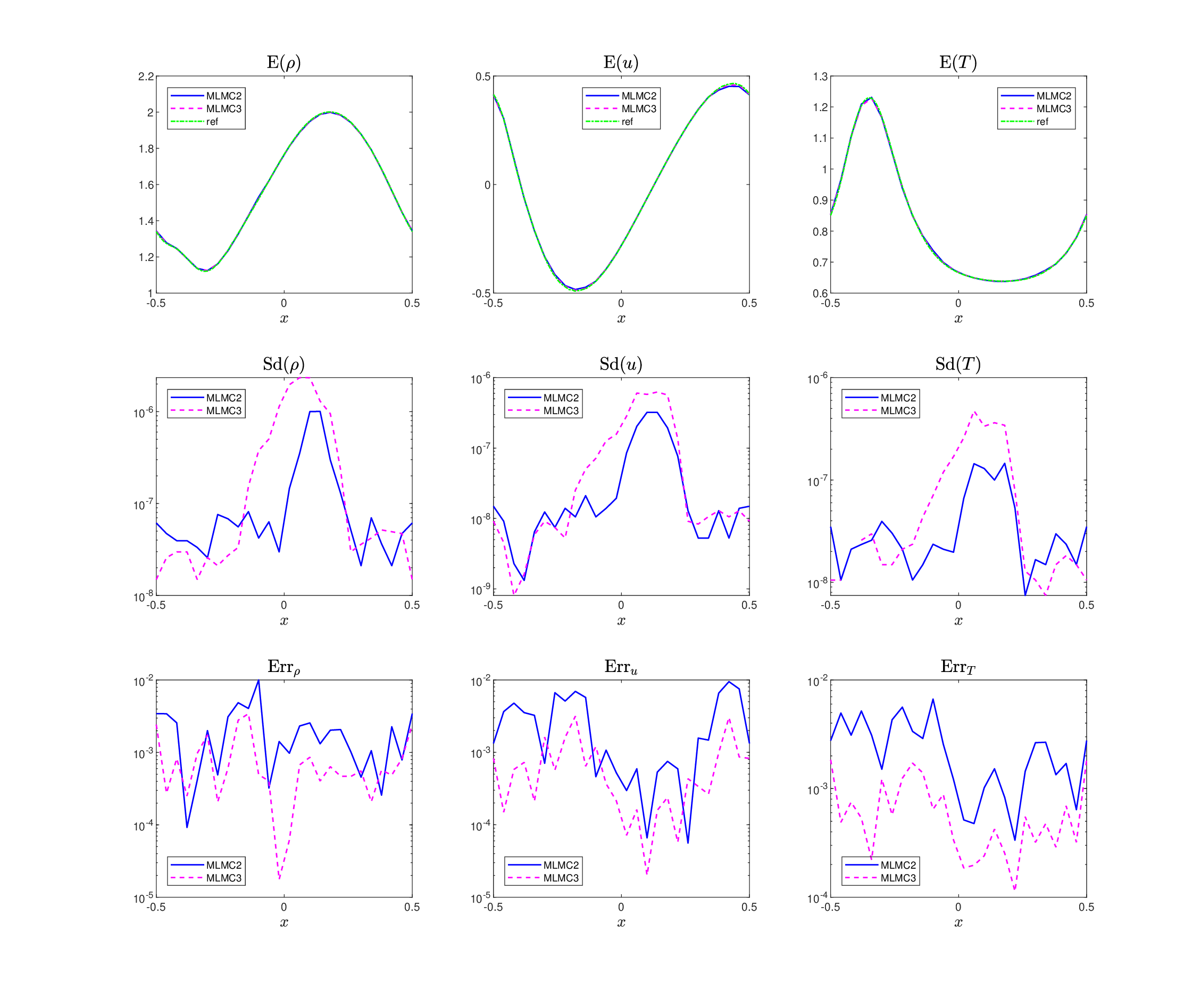}
	\caption{Test Four (b): Mixed regime with uncertainty in collision operator. Expectation, standard deviation and pointwise error of density $\rho$, mean velocity $u$ and temperature $T$ at time $t=0.2$ by the two-level MLMC ($N_0=50, N_1=100$) and the three-level MLMC ($N_0=50, N_1=100, N_2=200$).}\label{Fig:Ex5_Q_MLMC}
\end{figure}

Figure \ref{Fig:Ex5_Q_MLMC} compares the performance of the two-level ($N_0=25, N_1=50$) and the three-level MLMC ($N_0=25, N_1=50, N_2=100$) implementations by presenting the expectation, standard deviation, and pointwise error of density $\rho$, mean velocity $u$, and temperature $T$ at $t=0.2$. The results demonstrate that since the uncertainty is introduced solely in the collision operator, its impact is confined to kinetic regimes, while leaving the fluid regions unaffected. Therefore, the uncertainty propagation by the APH-based MLMC method remains localized to non-equilibrium regions; and this localized nature of collision-induced uncertainty leads to significantly smaller standard deviations compared to Test Four (a) where uncertainties are imposed in initial conditions.
Meanwhile, the three-level MLMC demonstrates superior accuracy over the two-level MLMC implementation, as clearly visible in the error plots, see the last row in Figure \ref{Fig:Ex5_Q_MLMC}. 
Therefore, the APH-based MLMC approach successfully distinguishes between uncertain kinetic regions and deterministic fluid zones, effectively handling the mixed-regime nature of the problem.

This test confirms the MLMC method's capability to efficiently resolve uncertainty in collision operators while maintaining accuracy through its multi-level structure, particularly beneficial for problems with localized stochastic effects. The results highlight how uncertainty propagation depends critically on its source (collision vs initial conditions), with collision uncertainty producing more localized effects that the APH-based MLMC method can exploit for computational efficiency.

\textbf{Test Four (c).} In this test case, we apply both the bi-fidelity and tri-fidelity approximations for the Boltzmann equation with uncertainty in initial conditions based on the hybrid method. 
For the bi-fidelity method, we select the Euler solver as the low-fidelity solver and the APH solver as the high-fidelity solver, instead of using the conventional AP solver for the Boltzmann equation as the high-fidelity solver. In the tri-fidelity method, we employ the Boltzmann solver as the high-fidelity solver, the APH solver as the medium-fidelity solver, and the Euler solver as the low-fidelity solver.

In the APH solver, the kinetic regimes are solved numerically with a spatial resolution of $\Delta x=0.01$, a temporal resolution of $\Delta t=8 \times 10^{-4}$, and $N_v=32$ velocity discretization points, up to the final time $t=0.2$. The fluid regimes in the APH solver employ the same spatial and temporal resolution as the kinetic regimes. 
For both the Boltzmann solver (high-fidelity model in the tri-fidelity method) and the Euler solver (the low-fidelity model in both bi- and tri-fidelity methods), the numerical parameters including temporal, spatial and velocity discretizations are identical to those used in the hybrid solver.

\begin{figure}[t]
	\centering
	\includegraphics[width=0.7\textwidth]{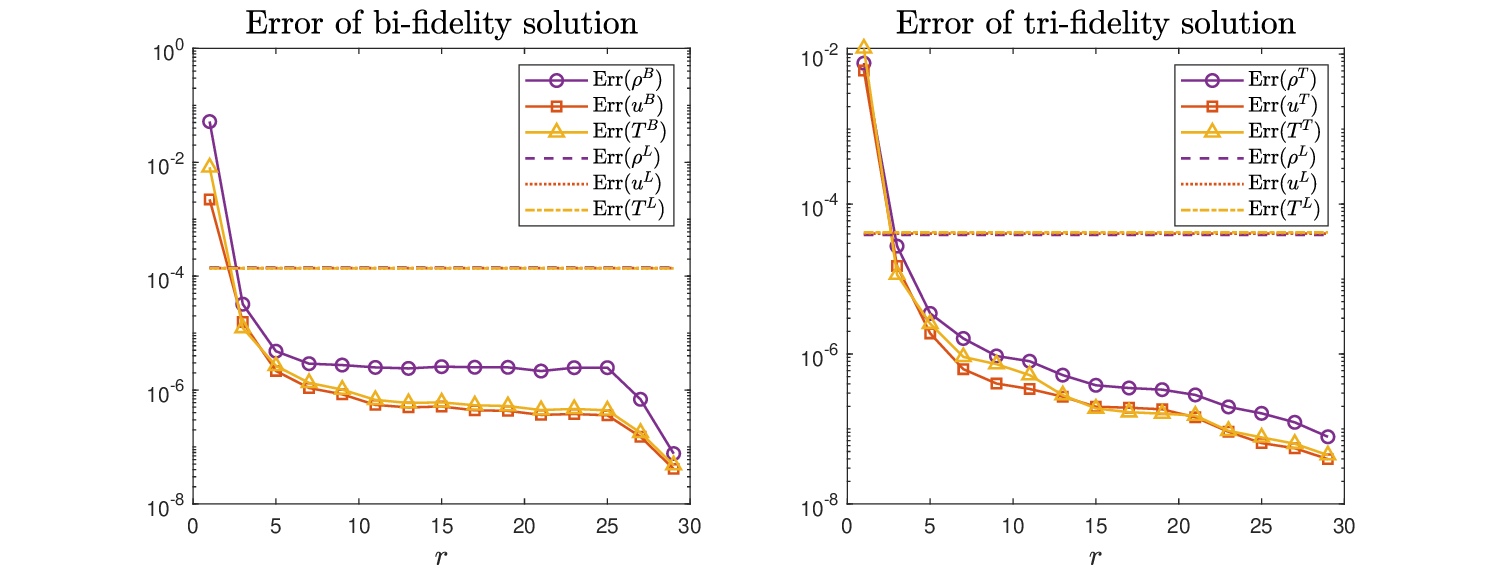}
	\caption{Test Four (c): mixed regime with uncertainty in initial conditions by multi-fidelity methods. The mean $L^2$ error of the bi- and tri-fidelity approximation of density $\rho$, mean velocity $u$ and temperature $T$ at times t = 0.2 with respect to the number of high-fidelity simulation runs.}\label{Fig:Ex5_bi_tri_err}
\end{figure}

The mean $L^2$ error between high-fidelity $\left(u^H\right)$ and multi-fidelity $u^\mathcal{F}$ solutions ($\mathcal{F}=B$ for bi-fidelity and $\mathcal{F}=T$ for tri-fidelity) at final time $t$ is evaluated by choosing a fixed set of points $\left\{\hat{z}_i\right\}_{i=1}^N \subset \Omega_z$ independent of the sampling sets $\Gamma$: 
$$
\mathcal{E}(u^H,u^\mathcal{F})= \frac{1}{N} \sum_{i=1}^N\left\|u^H\left(\hat{z}_i, t\right)-u^\mathcal{F}\left(\hat{z}_i, t\right)\right\|_{L^2(D)},
$$
where $\|\cdot\|_{L^2(D)}$ denotes the $L^2$-norm over the physical domain $D$.
Figure \ref{Fig:Ex5_bi_tri_err} plots the mean $L^2$ error of these approximations for the macroscopic quantities $\rho, u$, and $T$ with respect to the number of high-fidelity simulation runs. 
We observe a fast convergence of the $L^2$ errors. Notably, the dotted lines representing errors between high- and low-fidelity solutions are substantially larger than those between high- and multi-fidelity approximations. 
This demonstrates two key advantages:
\begin{itemize}
    \item [(1)] Bi-fidelity effectiveness: Using the APH solver as the high-fidelity model yields good bifidelity approximations. Therefore The APH solver provides a better high-fidelity choice for bi-fidelity approximation as \Cref{sec:deter} established the APH solver's superior efficiency over the full kinetic solver.
    \item [(2)] Tri-fidelity enhancement: Employing the APH solver as the medium-fidelity model provides dual benefits: it improves key sample point selection while enabling higher-accuracy tri-fidelity approximations compared to bi-fidelity approaches.
\end{itemize}

\begin{figure}[t]
	\centering
	\includegraphics[width=0.9\textwidth]{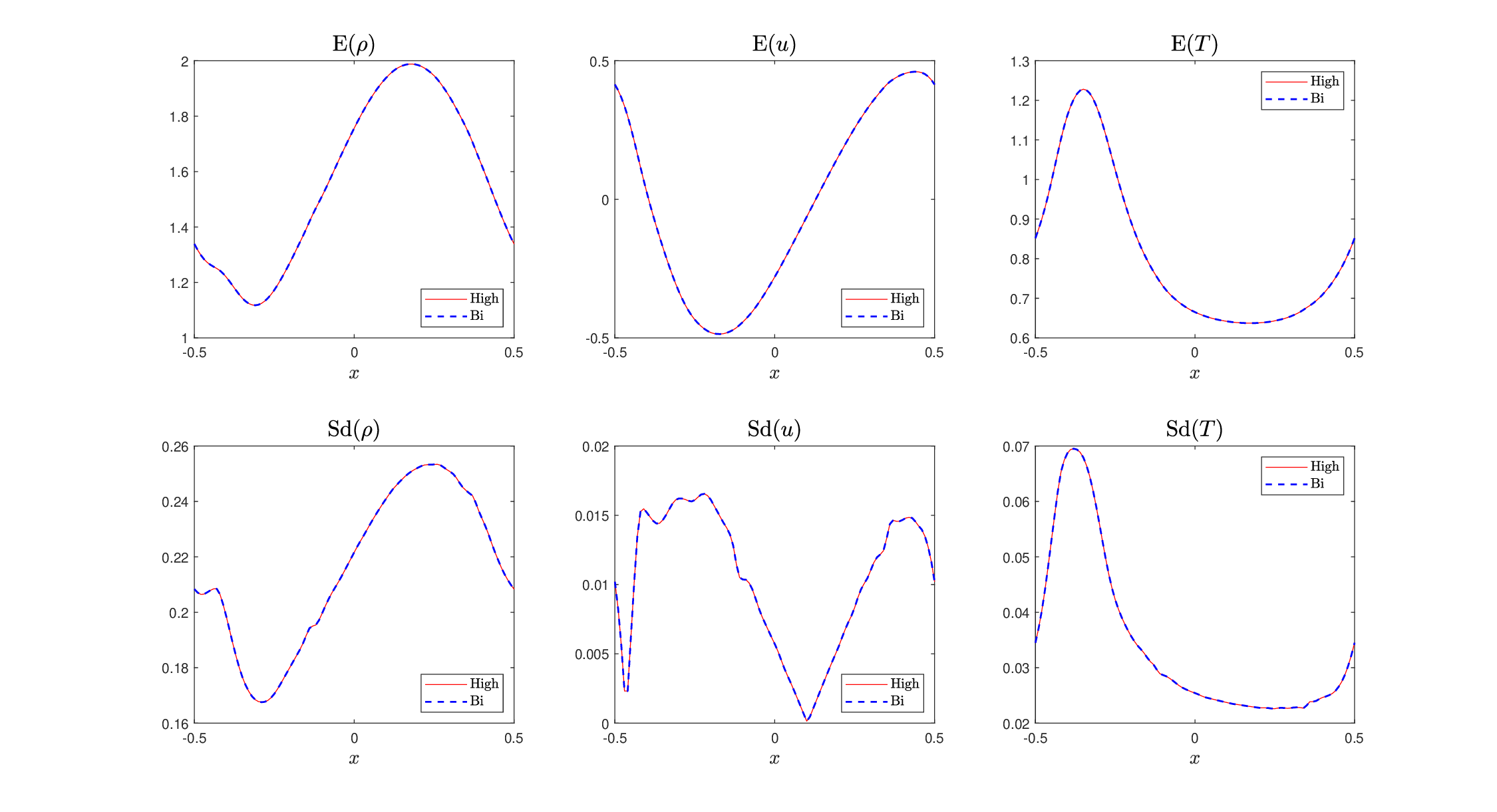}
	\includegraphics[width=0.9\textwidth]{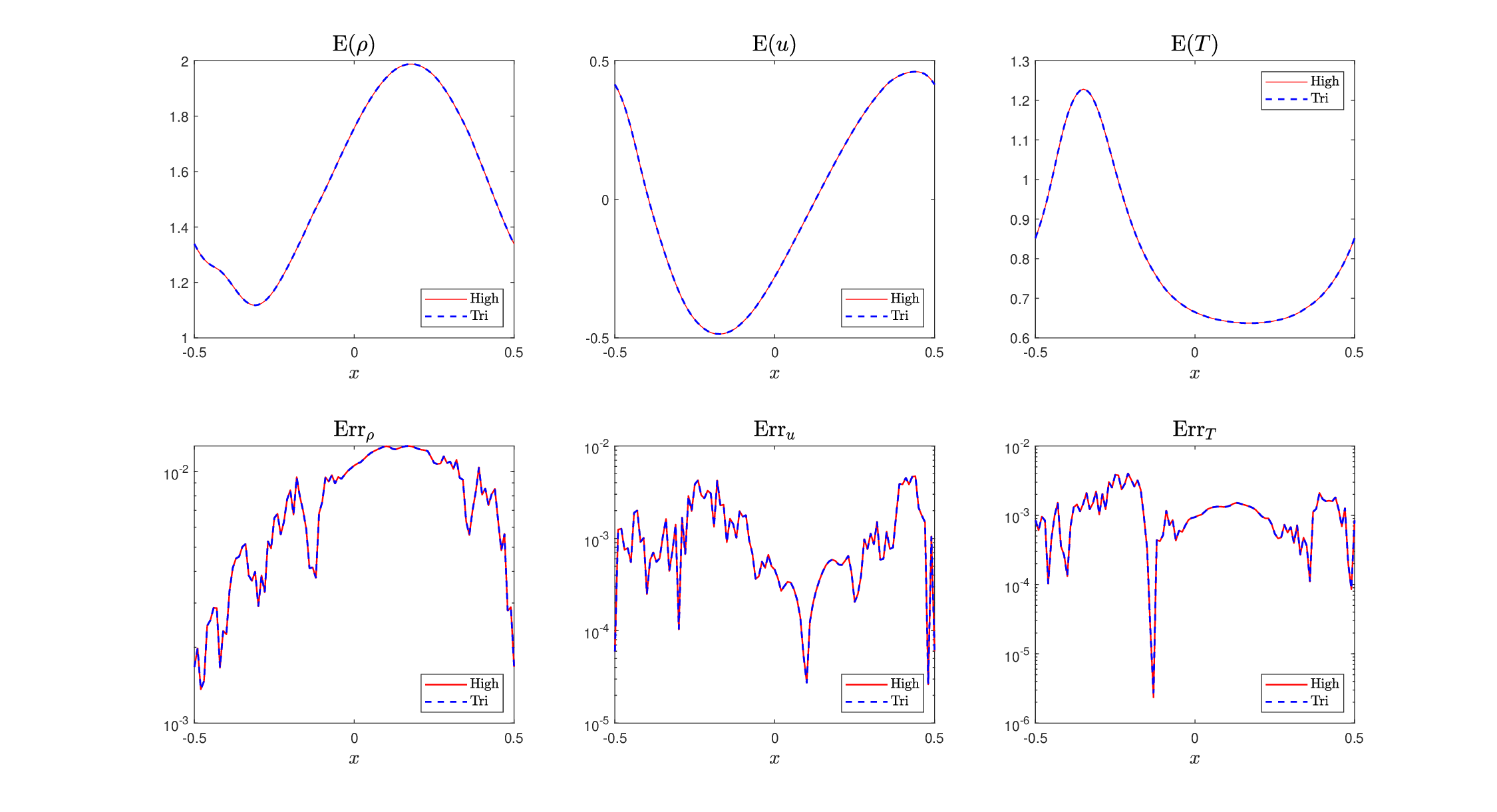}
	\caption{Test Four (c): mixed regime with uncertainty in initial conditions by multi-fidelity methods. The first two rows: Expectation and standard deviation of density $\rho$, mean velocity $u$ and temperature $T$ at times t = 0.2 by bi-fidelity approximation. The last two rows: Expectation and error of density $\rho$, mean velocity $u$ and temperature $T$ at times t = 0.2 by tri-fidelity approximation.}\label{Fig:Ex5_bi_tri_plot}
\end{figure}

In Figure \ref{Fig:Ex5_bi_tri_plot}, we present the expectation, standard deviation and pointwise error of the bi- and tri-fidelity approximations. The first two rows show the mean and standard deviation of $\rho, u$, and $T$ computed using 25 high-fidelity by the bi-fidelity method. The results show that with $N=1000$ low-fidelity runs of the Euler model, together with only 25 runs of the hybrid solver to the Boltzmann model, one can get the bi-fidelity solutions achieving excellent agreement with the high-fidelity solutions, which implies that the bi-fidelity approximation accurately captures the characteristics of the macroscopic quantities in the random space. However, relying solely on the low-fidelity Euler equations fails to maintain comparable accuracy, particularly across multiple scaling regimes where the Knudsen number $\varepsilon$ varies from $10^{-3}$ to 1. This limitation stems from the progressively deteriorating discrepancy between the Boltzmann solution and its Euler approximation as $\varepsilon$ increases. Such behavior underscores the crucial advantage of the bi-fidelity approach in bridging these multi-scale challenges.

Comparing Figure \ref{Fig:Ex5_bi_tri_plot} and Figure \ref{Fig:Ex5_MLMC3_100}, we observe that both the multi-fidelity method and the MLMC method can achieve comparable accuracy levels of $10^{-3}$. The computational complexity for the MLMC method is approximately $\mathcal{O}(25 \times M+50 \times M / 4$ $+100 \times M / 16 ) = \mathcal{O}(3 \times 25 \times M)$. For the multi-fidelity approach, the total work for the offline stage is $\mathcal{O}(100 \times M)$, while the online stage only needs a single low-fidelity computation to obtain the multi-fidelity approximation. These findings demonstrate that when the offline computation is permissible, the multi-fidelity method should be prioritized due to its computational efficiency in the online stage. However, when facing strict computational constraints, the control variate MLMC method serves as a nice alternative, providing a balance between accuracy and computational cost.
Notably, for sufficiently smooth solutions, MLMC can actually surpass bi-fidelity approaches in efficiency, achieving the accuracy of $10^{-3}$ with only 25 coarsest-level meshes and requiring just $M/16$ samples at 100-mesh resolution, while maintaining robust performance across various regimes.
This comparative analysis provides practical guidance for method selection based on problem characteristics and accuracy requirements.

\section{Conclusion}
\label{sec6:conclusion}

In this work, we numerically investigated the Boltzmann equation with uncertainties in the initial condition or collision operator using stochastic collocation-based uncertainty quantification (UQ) methods. A key component of our approach is a hierarchical hybrid method that efficiently handles the Boltzmann equation across different regimes. The hybrid method relies on two criteria: the first determines the transition from the fluid system to the kinetic equation based on macroscopic quantities, while the second governs the transition from the kinetic equation to its hydrodynamic limit by comparing the distribution function with its equilibrium state. 
To further enhance computational efficiency, we employ a control variate multilevel Monte Carlo (MLMC) method, which not only outperforms standard Monte Carlo but also maintains high accuracy in near-fluid regimes and discontinuous cases where conventional methods may fail. Additionally, we have developed a multi-fidelity collocation method, exemplified by the bi-fidelity approach, where the low-fidelity solution is derived from the compressible Euler equations while the high-fidelity solution is computed using the asymptotic preserving hybrid  (APH) solver for the Boltzmann equation. This approach accurately captures macroscopic variations at a computational cost far lower than that of the high-fidelity conventional asymptotic-preserving (AP) solver for the multiscale Boltzmann equation.
Despite these advancements, a rigorous error analysis for the APH-based MLMC and multi-fidelity methods remains an open question, as it depends on the parameter choices for the hybrid criteria and the solution’s smoothness, which are topics for future investigation.
Our framework is naturally extensible to other kinetic equations, combining deterministic phase-space discretizations with Monte Carlo sampling in the random space. 
For UQ in multiscale kinetic equations, the method selection should be guided by both solution regularity and computational constraints.

\section*{Acknowledgement}
The work of Y. Lin was partially supported by National Key R\&D Program of China (2020YFA0712000) and National Natural Science Foundation of China (12201404). L. Liu thanks the kind hospitality of Y. Lin during her visit to Shanghai Jiao Tong University. She acknowledges the support by National Key R\&D Program of China (2021YFA1001200), Ministry of Science and Technology in China, Early Career Scheme (24301021) and General Research Fund (14303022 \& 14301423) funded by Research Grants Council of Hong Kong. The numerical computations were run on the Siyuan-1 cluster supported by the Center for High Performance Computing at Shanghai Jiao Tong University. We thank Prof. Shi Jin for helpful discussion for the manuscript.

\bibliographystyle{plain} 
\bibliography{references}

\end{document}